\documentclass[psamsfonts]{amsart}
\usepackage{latexsym}
\usepackage{amscd, amsfonts, eucal, mathrsfs, amsmath, amssymb, amsthm}
\usepackage[OT2,OT1]{fontenc}

\input xy
\xyoption{all}

\newcommand{\field}[1]{\mathbf #1}
\newcommand{\mf}[1]{\mathfrak #1}

\newcommand{\ms}[1]{\mathscr #1}
\newcommand{\widebar}[1]{\overline{#1}}

\newcommand{\R}{\field R}

\newcommand{\Z}{\field Z}

\newcommand{\simto}{\stackrel{\sim}{\to}}

\renewcommand{\phi}{\varphi}

\renewcommand{\hom}{\operatorname{Hom}}
\newcommand{\shom}{\ms H\!om}

\newcommand{\rshom}{\mathbf{R}\shom}
\newcommand{\saut}{\ms A\!ut}

\DeclareMathOperator{\rhom}{\operatorname{{\bf R}Hom}}

\newcommand{\spec}{\operatorname{Spec}}
\newcommand{\spf}{\operatorname{Spf}}

\newcommand{\G}{\field G} 

\renewcommand{\H}{\operatorname{H}}

\DeclareMathOperator{\ext}{\operatorname{Ext}}
\newcommand{\sext}{\ms E\!xt}

\DeclareMathOperator{\sTor}{\ms T\!or}

\DeclareMathOperator{\coh}{\operatorname{Coh}}

\DeclareMathOperator{\D}{\operatorname{\bf D}}
\DeclareMathOperator{\ppD}{\operatorname{\bf D}_{p}^{\mathit b}}
\DeclareMathOperator{\ppugD}{\operatorname{\bf D}_{pug}^{\mathit b}}
\DeclareMathOperator{\sppugD}{\ms D_{\operatorname{pug}}^{\mathit b}}

\DeclareMathOperator*{\tensor}{\otimes}
\DeclareMathOperator*{\ltensor}{\stackrel{\field L}{\otimes}}
\newcommand{\im}{\operatorname{im}}

\newcommand{\surj}{\twoheadrightarrow}
\newcommand{\inj}{\hookrightarrow}
\newcommand{\id}{\operatorname{id}}

\DeclareMathOperator{\coker}{\operatorname{coker}}

\newcommand{\dirlim}{\varinjlim}
\newcommand{\invlim}{\varprojlim}

\DeclareMathOperator{\V}{\bf V}

\newtheorem{lem}{Lemma}[subsection]
\newtheorem{thm}[lem]{Theorem}
\newtheorem*{theorem}{Theorem}

\newtheorem{prop}[lem]{Proposition}

\newtheorem{cor}[lem]{Corollary}

\theoremstyle{definition}
\newtheorem{defn}[lem]{Definition}

\newtheorem{example}[lem]{Example}

\newtheorem{construction}[lem]{Construction}

\theoremstyle{remark}
\newtheorem{remark}[lem]{Remark}
\newtheorem{remarks}[lem]{Remarks}

\title{Moduli of complexes on a proper morphism}
\author{Max Lieblich}
\address{Department of Mathematics, Brown University, Box 1917, Providence RI 02912}
\curraddr{Department of Mathematics, Princeton University, Princeton NJ 08544}
\email{{\tt lieblich@math.princeton.edu}}
\thanks{Work on this paper was supported by a National Science Foundation Postdoctoral Fellowship.}
\date{}
\begin{document}

\bibliographystyle{plain}

\begin{abstract} Given a proper morphism $X\to S$, we show that a large class of objects in the derived category of $X$ naturally form an Artin stack locally of finite presentation over $S$.  This class includes $S$-flat coherent sheaves and, more generally, contains the collection of all $S$-flat objects which can appear in the heart of a reasonable sheaf of $t$-structures on $X$.  In this sense, this is the Mother of all Moduli Spaces (of sheaves).  The proof proceeds by studying the finite presentation properties, deformation theory, and Grothendieck existence theorem for objects in the derived category, and then applying Artin's representability theorem.
\end{abstract}

\maketitle

\tableofcontents


\section{Introduction}

Recent work (\cite{bondal-orlov}, \cite{bridgeland}, \cite{chen}) has indicated the usefulness 
of constructing moduli spaces of certain types of objects in the derived category 
$\D(X)$ of a variety $X$.  The goal of this paper is to provide 
general foundations for this theory by constructing an algebraic 
stack (in the sense of Artin) parametrizing ``all'' of the objects 
in $\D(X)$ which could possibly arise in geometry.  In particular, we will prove the 
following theorem (stated more precisely below as \ref{T:main 
theorem}), generalizing recent work of Inaba \cite{inaba}:

\begin{theorem} Let $X\to S$ be a proper flat morphism of finite 
presentation between algebraic spaces which is fppf-locally on $S$ representable by 
schemes.  The stack $\ms D$ of objects $E$ in $\D(X)$ which are relatively 
perfect over $S$ and such that $\ext^{i}(E_{s},E_{s})=0$ for all 
geometric points $s\to S$ and all $i<0$ is an algebraic stack locally 
of finite presentation over $S$.
\end{theorem}
Choosing an appropriate $t$-structure should morally define an algebraic substack tailored 
to the problem at hand.  Thus, in Bridgeland's situation \cite{bridgeland}, given a 
crepant resolution $X\to Z$ with $Z$ a Gorenstein terminal 
projective $3$-fold, Bridgeland's choice of $t$-structure yields  
a nice locally closed substack of $\ms D$ whose coarse moduli space 
gives a flop $Y\to Z$.  (Bridgeland also imposes a stability 
condition, which we will study in later work.)  We hope that this general mechanism will 
clarify the situation and enable us to prove more results along these 
lines.

Previous work on this problem was undertaken by Inaba \cite{inaba}.  
Given a flat projective morphism $X\to S$ of locally Noetherian schemes, 
he constructed an algebraic space locally of finite presentation $C\to S$ 
parametrizing ``simple complexes'' on $X/S$.  As we will show 
in section \ref{S:easy applications} below, there is an open substack 
$\ms C\subset\ms D$ which is naturally a $\G_{m}$-gerbe over his space 
$C$.  In fact, one can see that the Brauer class of this gerbe gives 
the obstruction to the existence of a universal object on his space.  
(Inaba's construction is slightly more general, in the sense that he 
requires only $\ext^{-1}(E_{s},E_{s})$ to vanish; this is what he 
gains by sheafifying the moduli problem.  This condition is rarely 
found in nature in the absence of the vanishing of all negative exts, 
so we have chosen not to treat it here.  The interested reader can 
check that our methods will also yield his result, with enough care.)

We assume that the 
structural morphism $X\to S$ is flat, as this is always satisfied in practice, and its absence would require the use of 
derived algebraic geometry (in particular, taking derived base 
changes over $S$).  Recent work of Lurie, To\"en-Vezzosi, 
Behrend, and others should extend 
the results of this paper to the derived context, where one can 
eliminate the flatness hypothesis.  Furthermore, derived methods hold the prospect of yielding more structure on the stack we construct even in the case of a flat structural morphism, e.g., a virtual structure sheaf (in one of the $\infty$-categories of sheaves of algebras which are used in the derived theory).  In particular, these methods would yield a natural approach to constructing a virtual fundamental class.  In fact, To\"en and Vaqui\'e have recently posted a preprint \cite{toen-vaquie} which (among other things) carries out this derived program for perfect complexes on a proper smooth scheme over a field.  Nevertheless, until derived methods have 
penetrated more deeply into the foundations of algebraic geometry, it seems worthwhile to develop our results using classical techniques.

\subsection{The structure of the proof}

Our proof of the theorem uses Artin's representability criterion \cite{artin}.  This involves three main steps: 1) checking that the stack of complexes in question is locally of finite presentation and is locally quasi-separated with separated diagonal, 2) understanding the infinitesimal deformation theory of complexes, and 3) studying the effectivity of formal deformations (Grothedieck's Existence Theorem).  After giving a precise definition of the complexes involved in our moduli problem, we take up 1) in section \ref{S:the moduli problem}.  In section \ref{S:deformation theory} we treat 2) and 3).  The discussion of the infinitesimal deformation theory proceeds by reducing to the affine case, where one can work explicitly with free resolutions; a touch of bootstrapping then yields the general case.  The effectivization of formal deformations works by first realizing any formal deformation in the derived category as a formal deformation of actual complexes and then using induction on the number of non-zero cohomology sheaves to reduce 3) to the corresponding statement for coherent sheaves.  Finally, we feed all of the parts into Artin's beautiful machine in section \ref{S:algebraicity}. 

\subsection{Acknowledgements}

I would like to thank Dan Abramovich for suggesting this problem,  
making me aware of the existing work of Inaba, and giving me numerous helpful 
comments on rough drafts of this paper.  I would also like to 
thank Jacob Lurie and Johan de Jong for many enlightening suggestions and discussions, and Bertrand To\"en for his comments.  Finally, I would like to thank the referee for carefully reading this paper and bringing to my attention an error in the original version.

\section{The moduli problem}\label{S:the moduli problem}

In this section, we describe the types of complexes which will be of 
interest to us.  Our goal is ultimately to yield Artin stacks 
parametrizing objects in the heart of a $t$-structure (satisfying 
sufficiently many hypotheses which have not yet been entirely 
understood).  By 
imposing various stability conditions, we can hope to produce very 
well-behaved stacks and algebraic spaces.

\subsection{Definitions}

To start, let $\mf X$ be any topos and $A\to B$ a map of (unital 
commutative) rings in $\mf X$.  (This amount of generality will not 
last long, so the reader can just as well imagine a scheme or 
algebraic space.  However, the reader well-versed in these matters 
will recognize the importance of making such a general definition when 
it comes to understanding what it means for there to be a universal 
family over the moduli stacks we construct below.)

\begin{defn} A complex $E\in\D(B)$ is \emph{$A$-perfect\/} if there is 
a covering family of the final object $\{U_{i}\to e_{\mf X}\}$ such that 
for each $i$, $E_{U_{i}}$ is 
quasi-isomorphic to a bounded complex of $A$-flat quasi-coherent 
$B$-modules of finite presentation.  If $E$ is globally on $X$ 
quasi-isomorphic to a bounded complex of $A$-flat quasi-coherent 
sheaves of finite presentation then $E$ is \emph{strictly 
$A$-perfect\/}.  When $A$ is understood, we will also use the term 
\emph{relatively perfect (over the base)\/}.
\end{defn}

The primordial occasion to study relatively perfect complexes is when given a 
morphism of ringed topoi $f:X\to S$; one takes $A=f^{-1}\ms O_{S}$ and 
$B=\ms O_{X}$.  This notion was first defined (slightly differently) in \cite{sga6}; the reader can check that for flat 
morphisms of finite type of locally Noetherian schemes $X\to S$, our 
notion of relatively perfect complexes agrees with that of Illusie and 
Grothendieck in \cite{sga6}.

When $\mf X$ is quasi-compact, it follows that any relatively perfect complex $E$ 
is bounded, but this is easily seen to be unnecessary in general.  
Any strictly relatively perfect complex is bounded, and any unbounded relatively 
perfect complex cannot be strictly relatively perfect.  
It is worth noting that relatively perfect complexes need not be 
perfect (even in 
the case of perfect morphisms $X\to S$).  For example, one can take $S$ to 
be a point and $X$ a variety over the point with a singular point.  
For $E$, one simply takes the structure sheaf of the singular point.  
If $X$ is smooth over $S$, however, one can easily see that any 
$S$-perfect complex is in fact perfect on $X$.  These facts will not 
be relevant for the purposes of this paper.

\begin{example} If $f:X\to S$ is a morphism of algebraic spaces, then we 
can take the (small or large) \'etale topos of $X$ for $\mf X$, $f^{-1}\ms O_{S}$ for $A$ 
and $\ms O_{X}$ for $B$.  In this case, a complex $E\in\D(X)$ is 
$A$-perfect if and only if it is \'etale-locally on $X$ 
quasi-isomorphic to a bounded complex of $S$-flat quasi-coherent 
$\ms O_{X}$-modules of finite presentation.  We will also call this notion 
\emph{$S$-perfect\/} in this case.  Of course, any bounded complex 
$P^{\bullet}$ consisting of coherent $S$-flat sheaves on $X$ is 
$S$-perfect.
\end{example}

We briefly recall three lemmas of Grothendieck which will be repeatedly useful in 
the sequel.

\begin{lem}\label{L:standard groth} 
Let $A$ be a local ring with residue field $k$ and $P^{-1}\to P^{0}\to P^{1}$ a 
complex of finite free $A$-modules such that $\H^{0}(P\tensor_{A} 
k)=0$.  There is a decomposition $P^{0}=I^{-1}\oplus K^{1}$ such 
that the map $P^{-1}\to P^{0}$ factors through a surjection 
$P^{-1}\surj I^{-1}$ and the map $P^{0}\to P^{1}$ factors through an 
injection $K^{1}\inj P^{1}$.  Both modules $I^{-1}$ and $K^{1}$ are 
finite free $A$-modules.
\end{lem}
\begin{proof} We give a proof by induction on the rank of $P^{0}$; we 
are not sure if this is the standard proof.  First, suppose 
$P^{0}\cong A$.  By Nakayama's Lemma, the only interesting case to 
consider is when the map $P^{0}\tensor k\to P^{1}\tensor k$ is 
injective; we wish to show that the map $P^{-1}\to P^{0}$ must be the 
zero map.  To see this, it suffices to show that any element of 
$P^{1}$ with non-trivial annihilator must lie in $\mf m_{A}P^{1}$.  
By projecting to the factors, we may assume that $P^{1}\cong A$.  But 
now the statement follows immediately from the fact that any element not in 
$\mf m_{A}$ is a unit!

To prove the general case, suppose first that $P^{0}\tensor k\to 
P^{1}\tensor k$ is injective.  It follows that the same holds for 
$P^{0}\to P^{1}$ by the usual criterion for injectivity involving the non-vanishing of 
a determinant (and the fact that $A$ is local).  So we may suppose that the rank of the image of 
$P^{-1}$ in $P^{0}\tensor k$ is non-zero.  
Let $v_{1},\ldots,v_{r}\in P^{-1}$ be 
elements whose images in $P^{0}\tensor k$ form a basis for the image 
of $P^{-1}\tensor k$.  We claim that the submodule of $P^{-1}$ 
generated by the $v_{i}$ is a summand $V$ mapping isomorphically onto 
a summand of $P^{0}$ in the kernel of $P^{0}\to P^{-1}$.  To see this, 
note that the projection $P^{0}\to P^{0}\tensor k\to\im(P^{-1}\tensor 
k)$ lifts to a map $P^{0}\to A^{r}$ such that the images of the 
$v_{i}$ map (by Nakayama's lemma) to a basis for $A^{r}$.  Now we can 
write the complex as
$$W\oplus V\to S\oplus V\to P^{1},$$
so that $W\to S\to P^{1}$ is another complex of finite free modules 
with vanishing $0$th cohomology over $k$.  Since the rank of $S$ is 
strictly smaller than the rank of $P^{0}$, we are done by induction.
\end{proof}

\begin{lem}\label{L:non-triv groth lem} Let $A$ be a ring, $X$ an $A$-scheme of finite 
presentation, and $P^{-1}\to P^{0}\to P^{1}$ a complex of finitely 
presented $A$-flat $\ms O_X$-modules.  Suppose $a\in\spec A$ is a point.  If 
$\H^{0}(P^{\bullet}\tensor_{A}\kappa(a))=0$ then there is an open subscheme $U\subset X$ containing $X\tensor_A\kappa(a)$  such that the complex 
$P^{0}/\im(P^{-1})|_U\to P^{1}|_U$ consists of finitely presented $A$-flat 
$\ms O_U$-modules, has universally vanishing $\H^{0}$, and computes $\tau_{\geq 0}P^{\bullet}|_U$ after all base changes.
\end{lem}
\begin{proof}[Sketch of proof] One reduces to the situation where $A$ 
and $X=\spec B$ are local and Noetherian.  In this case, one can prove by 
induction from the case of $A$-modules of finite length and the faithful flatness of completion along any radical ideal that 
$\H^{0}(P\tensor_{A}\widehat M)=0$ for any finite $A$-module $M$, the 
hat denoting completion; one then 
concludes that $\H^{0}(P\tensor M)=0$ for any $A$-module $M$.  Looking 
at the sequence $$P^{0}/\im P^{-1}\to P^{1}\to \coker(P^{0}\to 
P^{1})\to 0$$ and using the universal vanishing of cohomology and 
Krull's theorem, one concludes first that the left-hand arrow is injective, next that 
$P^{0}/\im P^{-1}$ is $A$-flat, and finally that its 
formation commutes with base change.
\end{proof}

\begin{lem}[cher \`a Cartan] Suppose $\phi:A\to B$ is a morphism of rings 
in a topos $\mf X$.  Given any $E\in\D(A)$ and $F\in\D(B)$, there is 
a natural isomorphism
$$\R\shom_{A}(E,F|^{\R}_{A})\simto\R\shom_{B}(E\ltensor_{A}B,F),$$
where $F|^{\R}_{A}$ denotes the derived restriction of scalars 
$\D(B)\to\D(A)$.
\end{lem}
\begin{proof} The reader who desires an explicit proof can proceed as 
follows: resolve $E$ by a $K$-flat complex of $A$-modules $K\to E$ and 
$F$ by a $K$-injective complex of $B$-modules $F\to I$.  It is easy 
to see (using the tools of \cite{spaltenstein}, for example) that 
the complex $\shom^{\bullet}(K,I)$ computes both sides of the equality.
\end{proof}
In particular, given a scheme $X$, a closed subscheme $\iota:Y\subset 
X$, an $\ms O_{X}$-module $M$, and an $\ms O_{Y}$-module $N$, one 
has $\ext^{i}_{X}(M,N)=\ext^{i}_{Y}({\field L}\iota^{\ast}M,N)$.  This will 
be used repeatedly with little or no comment below.

\begin{prop} Let $S=\spec A$ be an affine scheme and $X\to S$ a finitely 
presented flat morphism.  If $X$ posseses an ample family of invertible sheaves 
then any $S$-perfect complex $E$ on $X$ is strictly $S$-perfect.
\end{prop}
\begin{proof} Since $X$ is quasi-compact, the complex $E$ is 
bounded, so using the ample family of invertible sheaves we may represent $E$ by a bounded above complex 
$P^{\bullet}$ of locally free $\ms O_{X}$-modules.  We claim that a 
sufficiently negative truncation $\tau_{\geq n}P^{\bullet}$ represents 
$E$ by a bounded complex of $S$-flat quasi-coherent $\ms O_{X}$-modules 
of finite presentation.    
Since $E$ is $S$-perfect and bounded, there exist $a$ and $b$ such that for all $x\in\spec A$, one 
has $E_{x}\in\D^{[a,b]}(X_{x})$.  Applying \ref{L:non-triv groth lem}, 
one sees that $\tau_{\geq a}P^{\bullet}$ is the desired complex.
\end{proof}

\begin{cor}\label{C:strict on affine} If $X/S$ is a flat finitely presented quasi-projective morphism to 
an affine scheme, then any $S$-perfect complex on $X$ is strictly 
$S$-perfect.
\end{cor}
This fact enables Inaba to define his functor in \cite{inaba} using only strictly 
$S$-perfect complexes.  However, as we wish to make clear in this 
paper, Inaba's condition is ``morally'' a local condition, and all of 
the relevant facts about such complexes may be derived from the local 
properties.  

We will now focus on the case 
of a morphism of algebraic spaces.  Much of what follows can be 
significantly generalized, but we believe that one sacrifices a great 
deal of clarity for only marginal mathematical gain.

\begin{defn} Let $f:X\to S$ be a flat 
morphism of algebraic spaces.  An $S$-perfect complex $E\in\D(\ms O_{X})$ 
is \emph{gluable\/} if $\R f_{\ast}\R\shom(E,E)\in\D(\ms O_{S})^{\geq 
0}$.  It is \emph{universally gluable\/} if this remains true upon 
arbitrary base change $T\to S$.
\end{defn}
Note that universal gluability is equivalent to the vanishing of all 
negative cohomology for the ``global'' complex $\rhom(E_{T},E_{T})$ 
for any affine $T\to S$.  
Not surprisingly, given sufficiently many finiteness conditions there is a fiberwise criterion for universal 
gluability.

\begin{prop} Let $f:X\to S$ be a proper flat morphism of finite 
presentation between locally Noetherian algebraic spaces.  An $S$-perfect complex $E\in\D(X)$ is 
universally gluable if and only if $\ext^{-i}(E_{s},E_{s})=0$ for all 
geometric points $s\to S$ and all $i<0$. 
\end{prop}
\begin{proof} First suppose $S$ is the spectrum of a complete local 
Noetherian ring $(A,\mf m,k)$.  Given any triangle $M\to N\to P\to $ 
in $\D(A)$, there results a triangle $E\ltensor M\to E\ltensor N\to 
E\ltensor P\to $, whence one finds a triangle

$$\xymatrix@C=5pt{ &\R\shom(E,E\ltensor P)\ar[dl]_+& \\
\R f_{\ast}\R\shom(E,E\ltensor M)\ar[rr] & &\R f_{\ast}\R\shom(E,E\ltensor 
N).\ar[ul]}$$
Let $\iota:X_{k}\inj X$ be the natural closed immersion.  Since $$\R 
f_{\ast}\R\shom(E,E\ltensor k)=\R 
f_{\ast}\R\iota_{\ast}\R\shom_{X_{k}}(E_{k},E_{k}),$$ we deduce from 
the above triangle that for any $A$-module $M$ of finite length, we 
have $$(\ast)\quad \ext^{i}(E,E\ltensor M)=0$$ for all $i<0$.  Applying 
\ref{L:constructibility}(2) below, we see that $(\ast)$ holds for any finite 
$A$-module $M$.  In particular, it holds for $M=A$.  By faithfully 
flat descent, $(\ast)$ therefore holds for the sections over any 
localization of $S$; similar reasoning shows that it holds over any 
localization of any $T\to S$.  The result follows.
\end{proof}

\begin{prop} Given a flat morphism $f:X\to S$ of algebraic spaces 
as above, the fibered category of universally gluable $S$-perfect complexes on $X$ 
forms a stack on $S$ in the fpqc topology.
\end{prop}
\begin{proof} See Corollaire 2.1.23 \cite{BBD} or Theorem 2.1.9 of \cite{abramovich-polishchuk}.
\end{proof}

Given $f:X\to S$ as above, we will write $\ppD(X/S)$ for the category of 
$S$-perfect complexes on $X$, $\ppugD(X/S)$ for the category of 
universally gluable objects, and $\sppugD(X/S)$ for the stack of 
universally gluable $S$-perfect complexes on $X$.

\subsection{Finiteness and separation properties}

For unfortunate technical reasons, we assume in this section 
that $X\to S$ is fppf-locally representable by quasi-compact separated 
schemes.  This includes any algebraic space $X$ which is a flat form 
of a quasi-compact separated scheme over another scheme $S$.  Such spaces arise 
naturally e.g.\ when studying relative curves of genus $1$ without marked 
points.

\begin{prop}\label{P:local finite presentation} Let $X/R$ be a flat 
scheme of finite presentation with affine diagonal.  Let $R=\dirlim R_{\alpha}$ be a directed colimit of 
$S$-rings, and suppose $E\in\ppD(X_{R}/R)$.  There exists $\alpha$ and 
$E_{\alpha}\in\ppD(X_{R_{\alpha}})$ such that 
$E_{\alpha}\ltensor_{R_{\alpha}}R\cong E$.  Furthermore, given $E$ 
and $F$ in $\ppD(X_{R_{\alpha}}/R_{\alpha})$ and an isomorphism $\phi:E\simto F$ in 
$\ppD(X_{R}/R)$, there is some $\beta>\alpha$ and a morphism 
$\phi_{\beta}:E_{\beta}\simto F_{\beta}$ which pulls back to $\phi$. 
\end{prop}
\begin{proof} By the techniques of \S 8 of \cite{ega4-3} (especially Th\'eor\`eme 8.10.5), it is not difficult to see that there is an element $\alpha=0$ in the system indexing the colimit and an $X_0/R_0$ such that $X_0\tensor_{R_0}R\cong X$.  Furthermore, we are free to replace the system $\{\alpha\}$ by the cofinal system $\{\alpha\geq 0\}$.  We proceed by induction on the number of 
affines in an open cover of $X_{0}$.  To this end, write $X_{0}=U_{0}\cup 
V_{0}$ with $U_{0}$ affine and $V_{0}$ a union of fewer affines.  Let $i^{U_{0}}:U_{0}\inj X_{0}$ and 
$i^{V_{0}}:V_{0}\inj X_{0}$ be the natural inclusions.  There is a triangle
$$\R i^{U\cap V}_{!}E_{U\cap V}\to\R i^{U}_{!}E_{U}\oplus\R i^{V}_{!}E_{V}\to E\to.$$  
Furthermore, one can check that 
the exhibited triangle is compatible with base extension.
Thus, it suffices to show 
that $E_{U}$ and $E_{V}$ are defined over some $R_{\alpha}$ and that 
the map $\R i^{U\cap V}_{!}E_{U\cap V}\to\R i^{U}_{!}E_{U}\oplus\R 
i^{V}_{!}E_{V}$, which is just the adjoint of the natural map $E_{U\cap V}\to 
E_{U}|_{U\cap V}\oplus E_{V}|_{U\cap V}$, is defined over some 
$R_{\alpha}$.  By \ref{C:strict on affine}, on $U$ we may resolve $E$ 
by a bounded complex of $R$-flat quasi-coherent sheaves of finite 
presentation.  Thus, by standard limiting results of 
Grothendieck \cite{ega4-3}, 
there is some $\alpha$ such that $E_{U}$ is the derived 
extension of scalars of some 
$E_{U,\alpha}\in\ppD(X\tensor_{R_{0}}R_{\alpha})$.  By the induction 
hypothesis, the same is true for $E_{V}$.  

A proof that the map $E_{U\cap V}\to E_{U}|_{U\cap V}\oplus E_{V}|_{U\cap V}$
descends proceeds again by induction on the number of affines in a 
cover, using the same triangle as above.  In fact, this is 
subsumed in the second statement of the proposition, which we now prove.  
Suppose given $\phi:E\to F$; this is clearly compatible with the formation of the ``covering 
triangles'' shown above.  Consider the restrictions $E_{U}$ and 
$F_{U}$ to the affine $U$.  Representing $E$ and $F$ by truncations 
$\tau_{\geq n}P^{\bullet}$ and $\tau_{\geq n}Q^{\bullet}$ 
of bounded above complexes of finite free $\ms O_{U}$-modules, one 
easily sees that any morphism $E_{U}\to F_{U}$ is in fact represented 
by a map $\tau_{\geq n}P^{\bullet}\to\tau_{\geq n}Q^{\bullet}$.  
(Indeed, it is certainly represented by a map 
$P^{\bullet}\to\tau_{\geq n}Q^{\bullet}$, and this factors through the 
truncation.)  By the standard limiting results of Grothendieck, this 
comes from some finite stage.  Thus, the first two vertical arrows in 
the diagram of triangles
$$\xymatrix{\R i^{U\cap V}_{!}E_{U\cap V}\ar[r]\ar[d] & \R 
i^{U}_{!}E_{U}\oplus\R i^{V}_{!}E_{V}\ar[r]\ar[d] & E\ar[r]\ar[d] &\\
\R i^{U\cap V}_{!}F_{U\cap V}\ar[r] & \R 
i^{U}_{!}F_{U}\oplus\R i^{V}_{!}F_{V}\ar[r] & F\ar[r] &
}$$
may be assumed to come from some finite stage $\beta$.  By the axioms 
for triangulated categories, there is some extension $\tilde\phi:E_{\beta}\to 
F_{\beta}$ fitting into the diagram over $R_{\beta}$.  Furthermore, 
it is easy to see that any two arrows which complete the diagram 
differ by an element of $\hom(\R i^{U\cap V}_{!}E_{U\cap V}[1],\R 
i^{U}_{!}F_{U}\oplus\R i^{V}_{!}F_{V})$.  By adjunction and induction on the number 
of affines in a covering (which uses the fact that the diagonal of $X$ is affine!), 
any such arrow arises at a finite stage, completing the proof.
\end{proof}

\begin{cor} Let $X\to S$ be a proper flat morphism of finite 
presentation between algebraic spaces which is fppf-locally on $S$ 
representable by schemes.  Then the stack $\sppugD(X/S)$ is locally 
of finite presentation over $S$.
\end{cor}
\begin{proof} We may assume that $S=\spec R$ and that $R=\dirlim 
R_{\alpha}$, as above.  We may further assume that there is a 
minimal $\alpha$, say $\alpha=0$, and an fppf morphism $R_{0}\to 
R'_{0}$ such that 1) there is $X_{0}/R_{0}$ proper of finite 
presentation and flat such that 
$X=X_{0}\tensor_{R_{0}}R$, and 2) $X_{0}\tensor_{R_{0}}R'_{0}$ is a scheme.  Let 
$R_{\alpha}'=R_{\alpha}\tensor_{R_{0}}R_{0}'$.  Let $E\in\ppugD(X/R)$.  
By \ref{P:local finite presentation} there is some index $\beta$ and 
a complex $E'_{\beta}$ such that $E':=E\ltensor_{R}R'$ is 
$E'_{\beta}\ltensor_{R'_{\beta}}R'$.  Since $E\in\ppugD$, we can 
identify $E$ with $(E',\phi)$, where $\phi$ is the gluing datum for 
$E$ with respect to the fppf covering $\spec R'\to\spec R$.  
Applying \ref{P:local finite presentation} once more, we can 
descend $\phi$ to a finite level, and we can ensure that $\phi$ 
satisfies the cocycle condition at a finite level.  The result follows 
by descent of objects in $\ppugD$.
\end{proof}

Using the finite presentation result just proved, we can prove that 
the stack $\sppugD(X/S)$ is locally quasi-separated with separated 
diagonal (a term which seems not to have a name, but which is 
now essentially required of any algebraic stack \cite{l-mb}).

\begin{prop}\label{P:locally quasi-sep essence} Let $f:X\to S$ be a flat proper morphism of finite 
presentation to a 
regular algebraic space of everywhere bounded 
dimension (i.e., there exists $n$ such that for every point $s\to S$, 
the dimension of $\ms O_{S,s}^{sh}$ is at most $n$).  Suppose 
$E,F\in\ppD(X)$ and for every geometric point $s\to S$, one has 
$\ext^{-1}_{X_{s}}(E_{s},F_{s})=0$.  
\begin{enumerate} 
    \item There is a dense open 
$S^{0}\subset S$ over which the functor 
$H:\coh_{S^{0}}\to \coh_{S^{0}}$ sending $M$ to $$\sext^{0}(f; E,F\ltensor{\field L} 
f^{\ast}M):=\H^{0}(\R 
f_{\ast}\R\shom(E,F\ltensor{\field L} f^{\ast}M))$$
has the form $H(M)=\shom(Q,M)$ for some locally free coherent sheaf $Q$ 
on $S^{0}$.  Furthermore, the functor $$(T\to S^{0})\mapsto\sext^{0}(f_{T}; E_{T},F_{T})$$ 
is representable by $\V(Q)$.  Finally, the similar functor $\sext^{-1}(f; 
E,F)$ vanishes.

    \item   If $S$ is the spectrum of a discrete valuation ring, 
then the global sections $\ext^{0}(E,F)$ form a locally free 
$\ms O_{S}$-module.
\end{enumerate}
\end{prop}
\begin{proof} We may work \'etale locally on $S$ and assume that $S$ 
is the spectrum of a regular Noetherian ring $R$.  Furthermore, since 
we are working with relative Ext sheaves, we may restrict our 
attention to schemes $T\to S$ which are affine.  By assumption, any 
$\ms O_{S}$-module has local homological dimension at any point bounded 
above by $\dim R<\infty$.  We will repeatedly use this fact in the 
sequel; for the purposes of abbreviation, we call this $(\ast)$.  Since $F$ 
is bounded and $E$ is bounded above (and locally isomorphic to a 
bounded above complex of locally free modules), we have (using 
$(\ast)$) that 
$\rshom(E,F)$ is compatible with base change in $S$ and for all 
$\ms O_{S}$-modules $M$, we have $\rshom(E,F\ltensor{\field L} 
f^{\ast}M)=\rshom(E,F)\ltensor{\field L} f^{\ast}M$.  (Without using $(\ast)$, 
this is no longer true, contrary to what seems to be asserted in the proof of 
Lemma 5.3.4 of \cite{abramovich-vistoli}.  
Using the techniques of this 
paper, one can reprove Lemma 5.3.4 of [\emph{loc.\ cit.}], at least 
in the cases of interest to the authors.)  
We are thus interested in the $0$th hypercohomology of the complex $\shom(E,F)\ltensor{\field L} f^{\ast}M$ 
and of the derived base changes $\shom(E,F)_{T}$.  Write $\ms 
C=\shom(E,F)$; this is a bounded below complex with coherent 
cohomology sheaves.  Using $(\ast)$, one can see that there is an 
equality $\R f_{\ast}\R\shom(E,E\ltensor{\field L} f^{\ast}M)=\R 
f_{\ast}\R\shom(E,E)\ltensor M$ in $\D(S)$ (the projection formula).  
If $d=\dim R$, we then see that for all $g:T\to S$, ${\field L} g^{\ast}(\tau_{\leq 
d}\ms C)\to{\field L} g^{\ast}\ms C$ is a quasi-isomorphism in degrees $\leq 
0$ and similarly for $\tau_{\leq d}\ms C\ltensor{\field L} f^{\ast}M\to\ms 
C\ltensor{\field L} f^{\ast}M$.  So to understand the hypercohomology of ${\field L} 
f^{\ast}\ms C$ and $\ms C\ltensor{\field L} f^{\ast} M$, we may replace $\ms C$ 
by a bounded complex (which has the additional property that it 
remains bounded upon all base changes, a property which we will not 
use).  By generic flatness and boundedness of $\tau_{\leq d}\ms C$, 
we may shrink $S$ to $S^{0}$ and assume that all cohomology sheaves are $S$-flat (and coherent).  It follows 
from the usual cohomology and base change arguments (\cite{hartshorne} 
or \cite{mumford}) that the formation of $\R f_{\ast}\tau_{\leq d}\ms 
C$ commutes with base change over $S^{0}$.  Using \ref{L:standard 
groth} finishes the proof.

When $S=\spec R$ is the spectrum of a discrete valuation ring with 
residue field $k$, we see  
from the projection formula that $\H^{-1}((\R f_{\ast}\tau_{\leq d}\ms 
C)\ltensor_{R}k)=0$.  Using \ref{L:standard groth} or III.5.3 of 
\cite{sga6}, we see that $\H^{0}(\R f_{\ast}\shom(E,F))$ is a torsion 
free finite $R$-module, hence is locally free as $R$ has dimension $1$.
\end{proof}

\begin{defn} We will say that \emph{$\sext^{-1}(E,E)$ vanishes in 
fibers\/} if for all $s\to S$ one has $\ext^{-1}_{X_{s}}(E_{s},E_{s})=0$.
\end{defn}

\begin{cor}\label{C:locally quasi-sep pre} Let $f:X\to S$ be as in 
\ref{P:locally quasi-sep essence}, and let $E\in\ppD(X)$ such 
that $\sext^{-1}(E,E)$ vanishes in fibers.  Given an automorphism 
$\phi:E\to E$ in $\ppD(X)$ such that $\phi_{s}=\id$ for a dense set of 
points $s\to S$, we have $\phi=\id$ on all of $S$.
\end{cor}
\begin{proof} Considering $\phi-\id$, we see that it suffices to show that any endomorphism of $E$ which 
is $0$ in a dense set of fibers vanishes identically on a dense open.  This in turn follows 
immediately from the fact that for some $S^{0}\subset S$, the sheaf 
$\sext^{0}(f;E,E)$ is representable by $\V(Q)$ for locally 
free $Q$.
\end{proof}

We can now prove the key result, which will show that the stacks we 
construct are locally quasi-separated.

\begin{cor}\label{C:local quasi-sep} Let $f:X\to S$ be a proper morphism of finite presentation 
of schemes with $S=\spec R$ affine and reduced.  Suppose $E\in\ppD(X)$ and 
$\sext^{-1}(E,E)$ vanishes in fibers.  Given an automorphism of $E$ 
which is equal to the identity in a dense set of fibers of $f$, there 
is an open subscheme $U\subset S$ such that $\phi|_{U}=\id$.
\end{cor}
\begin{proof} Write $R$ as a colimit of finite type $\Z$-algebras.  
By \ref{P:local finite presentation}, we may assume that $X$, $E$, and 
$\phi$ are defined over some such subalgebra $R_{0}$.  Since $\spec 
R\to\spec R_{0}$ is dominant, we see that a 
dense set of points in $\spec R$ maps to a dense set of points in 
$\spec R_{0}$.  
Since the vanishing of $\phi-\id$ is geometric in fibers, we see that 
it suffices to prove the statement over $R_{0}$.  Since $R_{0}$ is a 
finite-type reduced $\Z$-algebra, there is a dense open which is 
regular and has bounded dimension at any point.  Thus, we may assume 
that $R_{0}$ is a regular finite-type $\Z$-algebra (with bounded 
dimension at every point).  The result follows immediately from 
\ref{C:locally quasi-sep pre}.
\end{proof}

In particular, this applies to universally gluable $S$-perfect 
complexes.  As we will see below, this implies that the diagonal of 
$\sppugD(X/S)$ is of finite type.  

\section{Deformation theory of complexes}\label{S:deformation theory}

The deformation theory of complexes can be made quite explicit when 
working on a flat projective morphism, as in 
\cite{inaba}.  Our goal in this section is to develop the 
theory in much greater generality.  The reader will note, however, 
that our approach does not work in an arbitrary topos.  We leave it 
as a question to the reader whether or not this theory is an avatar of 
a much more general (and therefore elegant) theory.

Throughout this section, $X\to S$ will denote a flat morphism of 
finite presentation of (quasi-separated) algebraic spaces, $0\to I\to 
A\to A_{0}\to 0$ will denote a square-zero extension of $S$-rings, 
and $E_{0}\in\ppD(X_{A_{0}}/A_{0})$ will be a given $A_{0}$-perfect 
complex on $X$.  Let $\iota:X_{A_{0}}\inj X_{A}$ be the natural closed 
immersion.  We will systematically write $\bullet\ltensor_{A}A_{0}$ 
for the functor ${\field L}\iota^{\ast}$ and (sloppily) write nothing for 
$\R\iota_{\ast}$.  
Adjunction provides a map ${\field L}\iota^{\ast}\R\iota_{\ast}E_{0}\to E_{0}$ 
which we will thus write $E_{0}\ltensor_{A}A_{0}\to E_{0}$.  The 
(homotopy) kernel of this all-important map will be denoted $Q$; thus, 
there is a natural triangle $Q\to E_{0}\ltensor_{A}A_{0}\to E_{0}\to $ 
in $\D(X_{A_{0}})$.  We will use the notation $K(X)$ to denote the category of (cohomologically indexed) complexes of sheaves of $\ms O$-modules on $X$.

\subsection{Statement of the result}

We will prove the following result.  A deformation of $E_{0}$ to 
$X_{A}$ is a complex $E$ on $X_{A}$ along with an isomorphism 
$E\ltensor_{A}A_{0}\simto E_{0}$.  We will properly define these 
objects in section \ref{S:preparations}.  Recall that a \emph{pseudo-torsor\/} under a group $G$ in a topos (including the topos of sets) is an object $T$ with a $G$-action such that a section gives rise via the action to an isomorphism $G\simto T$.  There is no requirement that local sections exist (this makes $T$ a torsor), so e.g.\ a pseudo-torsor in the category of sets may be the empty set.

\begin{thm}\label{T:deformation theory} Given $X,A,A_{0},I,E_{0}$ as above.
    
    \begin{enumerate}
        \item There is an element 
        $\omega(E_{0})\in\ext^{2}_{X_{A_{0}}}(E_{0},E_{0}\ltensor_{A_{0}}I)$ which 
	vanishes if and only if there is a deformation of $E_{0}$ to $X_{A}$.
    
        \item The set of deformations of $E_{0}$ to $X_{A}$ is a 
        pseudo-torsor under 
        $$\ext^{1}_{X_{A_{0}}}(E_{0},E_{0}\ltensor_{A_{0}}I).$$
    
        \item Suppose $E_{0}$ is gluable.  Given a deformation $E$ of $E_{0}$ to $X_{A}$, the set 
        of infinitesimal automorphisms of $E$ is a torsor under 
        $$\ext^{0}_{X_{A_{0}}}(E_{0},E_{0}\ltensor_{A_{0}}I).$$
    \end{enumerate}
    
\end{thm}

For future reference, we note the following immediate consequence of 
the cher \`a Cartan isomorphisms and the associativity of the derived 
tensor product.

\begin{cor}\label{C:deformation corollary} Given $X,A,A_{0},I,E_{0}$ as 
above.  Suppose $J\supset I$ annihilates $I$, and let $\widebar 
A=A/J$, $\widebar E=E_{0}\ltensor_{A_{0}}\widebar A$.
    
    \begin{enumerate}
        \item There is an element 
        $\omega(E_{0})\in\ext^{2}_{X_{\widebar A}}(\widebar E,\widebar 
        E\ltensor_{\widebar A}I)$ which 
	vanishes if and only if there is a deformation of $E_{0}$ to $X_{A}$.
    
        \item The set of deformations of $E_{0}$ to $X_{A}$ is a 
        pseudo-torsor under 
        $$\ext^{1}_{X_{\widebar A}}(\widebar E,\widebar 
        E\ltensor_{\widebar A}I).$$
    
        \item Suppose $E_{0}$ is gluable.  Given a deformation $E$ of $E_{0}$ to $X_{A}$, the set 
        of infinitesimal automorphisms of $E$ is a torsor under 
        $\ext^{0}_{X_{\widebar A}}(\widebar E,\widebar 
        E\ltensor_{\widebar A}I)$.
    \end{enumerate}

\end{cor}
Thus, if $(A,\mf m,k)$ is a local ring and $I$ has the property that $\mf 
m I=0$, the deformation theory is governed (as expected) by the 
cohomology of the restriction of $E_{0}$ to the closed fiber.  This 
is familiar from the classical deformation theory of sheaves and will 
prove useful when we study the cases of the Grothendieck Existence 
Theorem for complexes which will be useful for us.

\subsection{Preparations}\label{S:preparations}

In this section we compare the two obvious notions of deformation.  
Then we recall some book-keeping results about certain special 
functorial $K$-flat resolutions of complexes of sheaves (implicit in 
work of Spaltenstein).  (Recall that a complex of sheaves $A$ is \emph{$K$-flat\/} if for every acyclic complex $B$, the complex $A\tensor B$ is acylic.)

\begin{defn}\label{D:deformation} A \emph{deformation\/} of $E_{0}$ to $X_{A}$ is an 
object $E\in\D(X_{A})$ along with an isomorphism 
$E\ltensor_{A}A_{0}\simto E_{0}$.  An isomorphism of deformations is 
an isomorphism $E\simto E'$ which respects the isomorphisms with $E_{0}$ on 
$X_{A_{0}}$.
\end{defn}

Given an object $E$ and a map $\phi:E\ltensor_{A}A_{0}\to 
E_{0}$, there is an induced map $\tilde\phi:E\ltensor_{A}I\to 
E_{0}\ltensor_{A_{0}}I$ in $\D(X_{A_{0}})$.  Indeed, there is a 
natural isomorphism 
$$E\ltensor_{A}I\simto E\ltensor_{A}A_{0}\ltensor_{A_{0}}I$$
in $\D(X_{A_{0}})$ arising from associativity of the derived tensor 
product.  The map $\phi$ yields
$$E\ltensor_{A}A_{0}\ltensor_{A_{0}}I\to E_{0}\ltensor_{A_{0}}I.$$
It follows that if $\phi$ is an isomorphism then so is $\tilde\phi$.

\begin{defn}\label{D:coho def} A \emph{cohomological deformation\/} of $E_{0}$ to 
$X_{A}$ is a triangle $$E_{0}\ltensor_{A_{0}}I\to E\to E_{0}\to $$
such that the natural map of triangles 
$$\xymatrix{E\ltensor_{A}I\ar[r]\ar[d]^{\tilde\phi} & 
E\ar[r]\ar[d]^{\id} & E\ltensor_{A}A_{0}\ar[r]\ar[d]^{\phi} &\\
E_{0}\ltensor_{A_{0}}I\ar[r] & E\ar[r] & E_{0}\ar[r] &}$$
is an isomorphism.
An isomorphism of cohomological deformations is a map of triangles which 
is the identity on the two ends.
\end{defn}

Given $E_{0}$, we thus have two groupoids: the category of 
deformations, which we will temporarily denote $D$, and the category 
of cohomological deformations, which we will temporarily denote $C$.  

\begin{lem} There is a natural equivalence of categories $D\to C$ with 
a natural quasi-inverse $C\to D$.
\end{lem}
\begin{proof} The functor $D\to C$ arises as 
follows: Suppose $(E,\phi)$ is an object of $D$.  Since $\phi:E\ltensor_{A}A_{0}\simto E_{0}$ is an isomorphism, so is 
$\tilde\phi:E\ltensor_{A}I\simto E_{0}\ltensor_{A_{0}}I$, and we see 
that the natural map of triangles is thus an isomorphism.  The 
inverse functor $C\to D$ arises by simply forgetting everything but 
the right-most vertical arrow.  We leave it to the reader to check 
that these functors define an equivalence.
\end{proof}

\begin{lem}\label{L:perfection persists} Suppose $X/S$ is of finite 
presentation.  If a complex 
$E_{0}$ is in $\ppD(X_{A_{0}}/A_{0})$ then any deformation of $E_{0}$ 
to $X_{A}$ is in $\ppD(X_{A}/A)$. 
\end{lem}
\begin{proof} This follows immediately from \ref{L:non-triv groth 
lem} (after reducing to the case of affine $X$).
\end{proof}

\begin{defn}\label{D:good complex} A \emph{good complex\/} is a complex $F\in K(X)$ 
    such that $F^{i}$ has the form $\bigoplus j_{!}\ms O_{U}$ 
    for \'etale morphisms $U\to X$ with $U$ affine.  Given $E\in 
    K(X)$, a \emph{good resolution\/} of $E$ is a quasi-isomorphism 
    $F\to E$ with $F$ a good complex.
\end{defn}
Note that a good resolution has the property that it is free on 
stalks.  The reader will easily check that any bounded-above good 
complex is $K$-flat.  It is not the case, however, that any good 
complex is $K$-flat; the canonical example is given near the end of 
the introduction of  
\cite{spaltenstein} and is due to Dold.  For our purposes, a slight 
modification of this example will be ultimately more instructive: the 
complex of $\Z/4\Z$-modules $$F:0\to\cdots\to 0\to\Z/4\Z\to\Z/4\Z\to\cdots,$$ 
where each map is multiplication by $2$ and the $\Z/4\Z$ terms start 
at the index $0$.  There is a ``right resolution'' $\Z/2\Z\to F$, but 
one can see that $F$ cannot be used to compute derived tensor products 
of $\Z/2\Z$ over $\Z/4\Z$.

\begin{lem}\label{L:good tor} Let $F$ be a complex of 
$\ms O_{X_{A_{0}}}$-modules.  Then there is a natural 
isomorphism $$\sTor_{1}^{\ms O_{X_{A}}}(F,\ms O_{X_{A_{0}}})\cong I\tensor 
F.$$
\end{lem}
\begin{proof} Simply tensor the sequence $0\to I\to A\to A_{0}\to 0$ 
(pulled back to $X$) with $F$ (and use flatness of $X/S$).
\end{proof}

\begin{lem}\label{L:natural good} There is a functor $\rho:K(X)\to K(X)$ 
along with a morphism $\rho\to\id$ such that 
\begin{enumerate}
    \item for all $E\in K(X)$, $\rho(E)\to 
E$ is a good resolution;

    \item given any map from a bounded above good complex $F\to E$, there is a natural 
    lift $F\to\rho(E)$ which is termwise split;
    
    \item given a square zero extension $0\to I\to A\to A_{0}\to 0$ 
    of $S$-rings and a complex $E\in\ K(X_{A_{0}})$, there is a natural 
    map $\rho_{X_{A}}(E)\to\rho_{X_{A_{0}}}(E)$ which is a termwise 
    surjection of sheaves.
    
\end{enumerate}
In particular, given a quasi-isomorphism $E\to E'$, there is an 
induced quasi-isomorphism $\rho(E)\to\rho(E')$.
\end{lem}
\begin{proof} To construct $\rho$, it suffices by the techniques of 
Spaltenstein \cite{spaltenstein} or the more elementary (functorial) homotopy colimit 
constructions of Neeman and B\"okstedt \cite{neeman-bokstedt} to prove it for $E\in\D^{\leq 0}(X)$.  Indeed, any complex $F$ is the homotopy colimit of its truncations $\tau_{\leq n} F\to F$, and this homotopy colimit can be explicitly computed as the mapping cone of the map $\bigoplus\tau_{\leq n} F\to\bigoplus\tau_{\leq n} F$ which on the summand $\tau_{\leq m} F$ maps into $\tau_{\leq m} F\oplus\tau_{\leq m+1} F$ by the identity and the natural map $\tau_{\leq m} F=\tau_{\leq m}\tau_{\leq m+1} F\to\tau_{\leq m+1} F$.  The most notable feature of this construction is that it produces a \emph{functorial\/} realization of any complex (up to functorial quasi-isomorphism) as a mapping cone of direct sums of \emph{functorial\/} bounded above complexes.  
Since the mapping cone of any map of good complexes is good and the direct sum of good complexes is good, we see that it suffices to produce the desired functorial good resolution in the bounded above case, whence upon temporarily reindexing we may assume the complex is bounded above by $0$.  (We have just described the approach of \cite{neeman-bokstedt}; Spaltenstein gives a more rigid method in \cite{spaltenstein} which yields the same reduction.)

We may thus suppose that $E^{i}=0$ for all 
$i>0$.  Let $F^{0}=\oplus_{U,\sigma}j_{!}\ms O_{U}$, where $U,\sigma$ 
runs over all \'etale affines $U\to X$ along with a section of 
$E^{0}$ over $U$.  Let $C^{1}(E)=E\times_{E^{0}}F^{0}$; there is a 
natural map of complexes $C^{1}(E)\to E=:C^{0}(E)$.  It is easy to 
check that this natural map is a quasi-isomorphism.  Define 
$C^{2}(E)$ by setting the degree 0 term equal to $F^{0}$ and the 
degree 1 term equal to $\bigoplus_{U',\tau} j_{!}\ms O_{U'}$, where 
$U',\tau$ runs over sections of $C^{1}(E)^{-1}$.  In this way, one 
inductively defines $\rho$.  It is clear from the construction that 
the first part of the lemma holds.

To verify the second part, let $F\to E$ be any map from a bounded above good complex; without loss of generality, we may assume $F^i=0$ for $i>0$.  To 
give a map from $j_{!}\ms O_{U}$ to $E^{0}$ is the same as giving a 
global section of $E^{0}|_{U}$.  Furthermore, to give a composite map 
$j^{V}_{!}\ms O_{V}\to j^{U}_{!}\ms O_{U}\to E^{0}$ is to give a map $V\to 
U$ over $X$ and a section of $E^{0}|_{U}$ whose pullback to 
$E^{0}|_{V}$ along this map is the section of $E^{0}|_{V}$ determined 
by the map $j^{V}_{!}\ms O_{V}\to E^{0}$.  This shows that the map $F\to 
E$ naturally lifts to a map $F\to C^{1}(E)$.  Continuing in this 
manner and proceeding by induction yields a sequence of maps $F\to 
C^{n}(E)$ which stabilize on the components of index $>-n$ (here $n$ 
runs through positive integers).  This yields the natural map 
$F\to\rho(E)$, as required.

The third part of the lemma follows from the construction and the fact that the \'etale 
topology is naturally invariant under infinitesimal deformations.  
Thus, for example, when constructing $C^{1}(E)$ any section of $E^{0}$ 
over an \'etale affine $U\to X_{A_{0}}$ is the restriction of a 
(unique!) section over the deformation $\tilde U\to X_{A}$ with 
the same underlying topological space.  The rest of the proof proceeds by induction.
\end{proof}

Note that if $X$ is affine, any free resolution $P\to E$ is also a 
good resolution (with every open $U$ equal to $X$!), hence if $P$ is bounded above there is a 
lift $P\to\rho(E)$ over $E$.

\begin{construction}\label{L:Q map} As at the beginning of this section, given $E_0\in K(X)$, let $Q$ denote the homotopy kernel of the adjunction map, so that there is an exact triangle $Q\to E_0\ltensor_A A_0\to E_0\to $ in $\D(X_{A_0})$.  We construct a natural map $\theta(E_{0}):Q\to E_{0}\ltensor_{A_{0}}I[1]$ 
in $\D^{b}(X_{A_{0}})$ which, as we show in a moment, governs the deformation theory of $E_0$.

We give (unfortunately) a somewhat ad hoc derivation.  
Applying \ref{L:natural good} yields a $K$-flat resolution 
$E_{0}^{\bullet}$ of $E_{0}$ in 
$\D(X_{A_{0}})$ and a $K$-flat resolution $E^{\bullet}$ of $E_{0}$ in 
$\D(X_{A})$ such that there is a termwise surjective quasi-isomorphism $\sigma:E^{\bullet}\to 
E_{0}^{\bullet}$.  Furthermore, 
since any good $K$-flat resolution of $E_{0}$ admits a (termwise 
split) quasi-isomorphism to a good $K$-flat resolution of a chosen 
$K$-injective resolution of $E_{0}$, it will follow that the 
map $\theta$ we construct is independent of the choice of resolutions.
Letting $\widetilde{Q}=\ker\sigma$, we find an exact 
sequence of complexes (using \ref{L:good tor}) 
$$0\to 
E_{0}^{\bullet}\tensor_{A_{0}}I\to \widetilde{Q}\tensor_{A}A_{0}\to 
\ker(E^{\bullet}\tensor_{A}A_{0}\to E_{0}^{\bullet})\to 0.$$
This defines the map $\theta(E_0)$.
\end{construction}

\begin{remark} In $\D(X_A)$, we can split the adjunction map $E_0\ltensor_{A} A_0\to E_0$ by tensoring the triangle $I\to A\to A_0\to$ with $E_0$ over $A$.  This yields an isomorphism $Q\to E_0\ltensor_A I[1]=E_0\ltensor_{A_0}I\ltensor_A A_0[1]$.  The reader will note that the image of $\theta$  
in $\D(X_{A})$ is the image of this natural $A$-linear isomorphism 
$Q\to E_{0}\ltensor_{A_{0}}I\ltensor_{A}A_{0}[1]$ under the natural 
adjunction $$E_{0}\ltensor_{A_{0}}I\ltensor_{A}A_{0}[1]\to 
E_{0}\ltensor_{A_{0}}I[1].$$  However, since $\D(X_{A_{0}})\to\D(X_{A})$ is not faithful, this is not sufficient to characterize $\theta(E_0)$.
\end{remark}

Ultimately, theorem \ref{T:deformation 
theory} will come about by applying 
$\rhom(\,\bullet\,,E_{0}\ltensor_{A_{0}}I)$ to the triangle $Q\to 
E_{0}\ltensor_{A}A_{0}\to E_{0}\to$: the element $\omega(E_{0})$ is 
the image of the element $\theta(E_{0})$ of \ref{L:Q map} under the 
coboundary, and the space of deformations is the fiber of the 
coboundary over $\theta$.  Proving this, 
however, will take a bit of trickery.

\begin{lem}\label{L:deformations live here} Let $E$ be a 
cohomological deformation of $E_{0}$, represented by 
$$\gamma\in\ext^{1}_{X_{A}}(E_{0},E_{0}\ltensor_{A_{0}}I)=
\hom_{X_{A_{0}}}(E_{0}\ltensor_{A}A_{0},E_{0}\ltensor_{A_{0}}I[1]).$$  
The image of $\gamma$ in $\hom_{X_{A_{0}}}(Q,E_{0}\ltensor_{A_{0}}I[1])$ 
via the map $Q\to E_{0}\ltensor_{A}A_{0}$ is equal to $\theta(E_{0})$.
\end{lem}
\begin{proof} Again, we give an ad hoc proof.  Let (by abuse of 
notation) $E_{0}$ be a good complex on $X_{A_{0}}$ representing 
$E_{0}$, $\widetilde E_{0}$ a good complex on $X_{A}$ representing 
$E_{0}$, and $E$ a good complex on $X_{A}$ representing a deformation 
of $E_{0}$ to $X_{A}$.  We may assume (passing to a mapping cylinder 
if necessary) that $\widetilde E_{0}\to E_{0}$ is surjective and 
$E\to \widetilde E_{0}$ is injective.  There results a diagram
$$\xymatrix{ & 0\ar[d] & 0\ar[d] & &\\
0\ar[r] & R\ar[r]\ar[d] & Q\ar[r]\ar[d] & K\ar[r]\ar[d] & 0\\
0\ar[r] & E\ar[r]\ar[d] & \widetilde E_{0}\ar[r]\ar[d] & K\ar[r] & 0\\
        & E_{0}\ar[r]\ar[d] & E_{0}\ar[d] & & \\
	& 0                 & 0,           & &}$$
where $K$ computes $E_{0}\ltensor_{A_{0}}I[1]$, $Q$ is acyclic, and $R$ 
computes $E_{0}\ltensor_{A_{0}}I$.  Furthermore, since $K$ is the 
mapping cylinder of a map of good complexes, we see that it is free on 
stalks.  Tensoring the diagram with $A_{0}$ yields
$$\xymatrix{  & 0\ar[d]        & 0\ar[d]        &          & \\
              & I\tensor E_{0}\ar[r]\ar[d]^{\id} & I\tensor E_{0}\ar[d] &            & \\
	    0\ar[r] & \widebar R\ar[r]\ar[d] & \widebar Q\ar[r]\ar[d] & \widebar 
	    K\ar[r]\ar[d]^{\id} & 0 \\
	    0\ar[r] & \widebar E\ar[r]\ar[d]&\widebar{\widetilde E}_{0}\ar[r]\ar[d] & 
	    \widebar K\ar[r] & 0 
	    \\
	    & E_{0}\ar[r]\ar[d] & E_{0}\ar[d] & &\\
	    & 0 & 0 & &}$$
Since $\widebar E\to E_{0}$ is a quasi-isomorphism, we conclude that 
$R\to\widebar R$ is a quasi-isomorphism, and therefore that the map in 
the derived category $\widebar K\to\widebar R[1]$ is identified with 
the adjunction map $E_{0}\ltensor_{A}I\to E_{0}\ltensor_{A_{0}}I$ via 
the natural isomorphisms $K\simto E_{0}\ltensor_{A_{0}}I$ and 
$I\tensor E_{0}\simto \widebar R$ induced by 
the isomorphism $\widebar E\simto E_{0}$ (of which the second is 
present in the diagram) and the requirement on the triangles arising 
in the definition of cohomological deformations \ref{D:coho def}.  Letting $\widetilde Q$ denote the kernel of 
the map $\widebar{\widetilde E}_{0}\to E_{0}$, one easily sees using 
the fact that $\widebar E\to E_{0}$ is a quasi-isomorphism that the 
induced map $\widetilde Q\to\widebar K$ is a quasi-isomorphism.  The 
equality follows from the diagram
$$\xymatrix{ & I\tensor E_{0}\ar[dl]^{\cong}\ar[d] &\\
\widebar R\ar[r] & \widebar Q\ar[d]\ar[r] & \widebar K\\
& \widetilde Q.\ar[ur]^{\cong} &}$$
\end{proof}

\begin{lem}\label{L:splitness} If there is a deformation 
$E_{0}\ltensor_{A_{0}}I\to E\to E_{0}\to $, then the triangle $Q\to 
E_{0}\ltensor_{A}A_{0}\to E_{0}\to$ is split.
\end{lem}
\begin{proof} The map $E\ltensor_{A}A_{0}\to E_{0}\ltensor_{A}A_{0}\to 
E_{0}$ gives a splitting of the adjunction map in $\D(X_{A_{0}})$.
\end{proof}
\begin{cor} When there is a deformation, the fiber of the 
natural (digaonal) map
$$\xymatrix{\hom_{X_{A}}(E_{0},E_{0}\ltensor_{A_{0}}I[1])\ar@{=}[d]\ar[dr]& \\
             \hom_{X_{A_{0}}}(E_{0}\ltensor_{A}A_{0},E_{0}\ltensor_{A_{0}}I[1])\ar[r]   & \hom_{X_{A_0}}(Q,E_{0}\ltensor_{A_{0}}I[1]) }$$
is a torsor under 
$\hom_{X_{A_{0}}}(E_{0},E_{0}\ltensor_{A_{0}}I[1])=\ext^{1}_{X_{A_{0}}}(E_{0},E_{0}\ltensor_{A_{0}}I)$
\end{cor}

\subsection{Complexes over an affine}

In this section, we assume $X$ is affine over $S$; in the deformation 
situation, we will have $X=\spec B$ 
with $B$ a flat $A$-algebra.  In this case, we 
can take resolutions of $E_{0}$ by bounded above complexes of free 
$B$-modules (or $B\tensor_{A}A_{0}$-modules).  This simplifies the 
picture significantly.

\begin{construction}\label{Constr:obstruction} Let $P_{0}$ be a bounded above complex of free 
$B_{0}$-modules representing $E_{0}$.  Lift each $P_{0}^{i}$ to a 
free $B$-module $P^{i}$.  Using projectivity, the differential 
$\bar d^{i}:P_{0}^{i}\to P_{0}^{i+1}$ may be lifted to some map of $B$-modules 
$d^{i}:P^{i}\to P^{i+1}$.  Since the $\bar d^{i}$ yield a complex 
$P_{0}$ and $I^{2}=0$, it is immediate that the maps $d^{i+1}\circ 
d^{i}:P^{i}\to P^{i+2}$ yield maps $\delta^{i}:P_{0}^{i}\to I\tensor 
P_{0}^{i+2}$, giving rise to a map of complexes $P_{0}\to I\tensor 
P_{0}[2]$.  The reader can check that the homotopy class of this map 
is independent of the choice of the lifts $d^{i}$.
Since $P_{0}$ is a resolution by free modules, this yields a unique 
element 
$$\omega(E_{0})\in\ext^{2}(P_{0},P_{0}\tensor_{A_{0}}I)=\ext^{2}(E_{0},E_{0}\ltensor_{A_{0}}I).$$
\end{construction}
We will verify in a moment that the element $\omega$ is independent 
of the choice of resolution $P_{0}\to E_{0}$.

\begin{lem}\label{L:affine obstruction} The image of $\theta(E_{0})$ under the coboundary 
$$\ext^{1}_{B_{0}}(Q,E_{0}\ltensor_{A_{0}}I)\to\ext^{2}_{B_{0}}(E_{0},E_{0}\ltensor_{A_{0}}I)$$
    is the element $\omega(E_{0})$ constructed above.
\end{lem}
\begin{proof} Choose a termwise surjective quasi-isomorphism $M\to P_{0}$ with $P_{0}$ a complex of free 
$B_{0}$-modules representing $E_{0}$ and $M$ a complex of 
free $B$-modules representing the image of $E_{0}$ in $\D(X)$.  The 
image of $\theta$ arises from the exact sequence associated to 
tensoring $0\to 
Q\to M\to P_{0}\to 0$ with $A_{0}$:
$$0\to I\tensor P_{0}\to \widebar Q\to\widebar M\to P_{0}\to 0.$$
We proceed to compute the realization of the map.  For any index 
$i$, let $P^{i}$ be the natural lift of $P_{0}^{i}$.  We can find maps 
$M^{i}\to P^{i}$ and $P^{i}\to M^{i}$ over $P_{0}^{i}$.  By 
Nakayama's Lemma, the composition $P^{i}\to M^{i}\to P^{i}$ is an 
isomorphism; thus, adjusting by an isomorphism, we may assume that 
the map $M^{i}\to P_{0}^{i}$ factors as a split surjection 
$M^{i}\surj P^{i}\to P_{0}^{i}$.  The splitting $P^{i}\to M^{i}$ 
gives rise to a splitting of $Q^{i}=\ker(M^{i}\to P_{0}^{i})$ as $I\tensor 
P^{i}\oplus N^{i}$, with locally free $N^{i}$.  There results for 
any $i$ a diagram
$$\xymatrix{0\ar[r] & I\tensor P_{0}^{i}\ar[r]\ar[d] & I\tensor 
P_{0}^{i}\oplus \widebar N^{i}\ar[r]\ar[d] & \widebar 
M^{i}\ar[r]\ar[d] & P_{0}^{i}\ar[r]\ar[d] & 0\\
0\ar[r] & I\tensor P_{0}^{i+1}\ar[r]\ar[d] & I\tensor 
P_{0}^{i+1}\oplus \widebar N^{i}\ar[r]\ar[d] & \widebar 
M^{i+1}\ar[r]\ar[d] & P_{0}^{i+1}\ar[r]\ar[d] & 0\\
0\ar[r] & I\tensor P_{0}^{i+2}\ar[r] & I\tensor 
P_{0}^{i+2}\oplus \widebar N^{i+2}\ar[r] & \widebar 
M^{i+2}\ar[r] & P_{0}^{i+2}\ar[r] & 0}$$
in which the direct sum splitting at each level is non-canonical and 
\emph{is not compatible with the coboundary\/}.  (As we will see in a 
moment, this lack of compatibility is precisely the point!)  
Furthermore, we have that $\widebar M^{i}$ splits (non-canonically) as 
$P_{0}^{i}\oplus\widebar N^{i}$.  Let $\sigma_{i}:P_{0}^{i}\to 
\widebar M^{i}$ be a choice of splitting; one could for example arrive 
at such a map by choosing a (non-linear) lift of each element into 
$P^{i}$, apply a splitting $P^{i}\to M^{i}$, and then take the image 
in $\widebar M^{i}$.  We will in fact assume that the splitting in question 
has this form.  To compute the map in the derived category under the 
hypothesis that the various complexes are termwise split is easy: 
given $t\in P_{0}^{i}$, first one forms 
$\sigma^{i+1}(d_{P_{0}}^{i}(t))-d_{\widebar M}^{i}(\sigma^{i}(t))$.  
This is an element $\beta(t)$ of $\widebar N^{i+1}$ by construction.  
Now one takes $d(\beta)\in I\tensor P_{0}^{i+1}\oplus \widebar 
N^{i+1}$; write $d(\beta)=f(t)\oplus g(t)$.  Setting the image of $t$ in 
$I\tensor P_{0}^{i+1}$ to be $f(t)$, there results a map $P_{0}\to 
I\tensor P_{0}[2]$, which 
realizes $\delta\theta$.  Another way to describe $f(t)$ is as 
follows: lift $t$ to an element of $P^{i}$.  Using the splittings 
$\tau_{i}:P^{i}\to M^{i}$ and the projections $\pi_{i}:M^{i}\to 
P^{i}$ (over $P_{0}$), define a map $d^{i}:P^{i}\to P^{i+1}$ by 
$\pi_{i+1}d_{M}^{i}\tau_{i}$.  One easily checks that $d^{i}$ is a 
map over the differential $P_{0}^{i}\to P_{0}^{i+1}$.  
We see immediately that 
$\tau_{i+1}d^{i}(t)-d_{M}^{i}(\tau_{i}(t))=\gamma(t)\in N^{i+1}$.  Applying 
$d^{i+1}_{M}$  to $\gamma(t)$ and projecting to the $P^{i+2}$ 
component yields the same result as applying $d^{i+1}$ to 
$d^{i}(t)$ (by the way we have defined the $d^{i}$!).  It follows 
that $\delta\theta=\omega$ and that consequently the construction 
of $\omega$ is independent of $P_{0}$.  (We already know that it is 
independent of the choice of lifts $d^{i}$, so we are free to take a 
convenient choice as we have done in this proof.)
\end{proof}

\begin{remark}\label{R:from projective to flat} The reader will note 
that the proof of \ref{L:affine obstruction} also works in a more 
general situation.  Given any $X$ and any $E_{0}$ (without 
boundedness conditions), by \ref{L:natural good} there is a termwise 
surjective quasi-isomorphism $M\to P_{0}$ with $P_{0}$ a good 
resolution of $E_{0}$ on $X_{A_{0}}$ and $M$ a good resolution of 
$E_{0}$ on $X_{A}$.  Furthermore, there is a natural lift $P^{i}$ of 
any term $P_{0}^{i}$ to a good sheaf on $X_{A}$.  By \ref{L:natural 
good}(2), there is a lift $P^{i}\to M^{i}$ over $P_{0}^{i}$ which is a 
split injection.  Thus, we find the same termwise splitting 
$P^{i}\inj M^{i}\surj P^{i}\to P_{0}^{i}$ as in \ref{L:affine 
obstruction}, and the argument produces a map $P_{0}\to P_{0}\tensor 
I[2]$ whose image in $\ext^{2}_{X_{A_{0}}}(E_{0},E_{0}\ltensor_{A_{0}}I)$ represents 
$\omega(E_{0})$.  This point will be useful in section \ref{S:groth 
exist} when we study the Grothendieck Existence theorem for relatively 
perfect complexes on $X/S$.
\end{remark}

\begin{prop}\label{P:affine deformation theory} Theorem \ref{T:deformation theory} holds for affine $X$.
\end{prop}
\begin{proof} From the construction \ref{Constr:obstruction}, it 
easily follows that if $\omega=0$ then there is a lift of the $\bar 
d^{i}$ to maps $d^{i}$ giving a structure of complex to the $P^{i}$; 
the resulting complex is easily seen to be a deformation of $E_{0}$.  
Furthermore, it is easy to see that the space of homotopy 
classes of lifts of the $\bar d^{i}$ to give such complexes with terms 
$P^{i}$ is principal homogeneous under $\ext^{1}(P_{0},P_{0}\tensor 
I)=\ext^{1}_{X_{A_{0}}}(E_{0},E_{0}\ltensor_{A_{0}}I).$  Using \ref{L:deformations 
live here}, it now follows that \emph{all\/} elements of 
$\ext^{1}_{X_{A}}(E_{0},E_{0}\ltensor_{A_{0}}I)$ mapping to $\theta$ 
arise from deformations.  The statement on infinitesimal 
automorphisms is left to the reader.
\end{proof}

\begin{cor}\label{C:explicit affine} Given a complex $P_{0}$ of free $B_{0}$-modules 
representing $E_{0}$, every deformation of $E_{0}$ to $B$ arises as a 
complex whose terms are the natural lifts $P^{i}$ of the $P_{0}^{i}$.
\end{cor}
\begin{proof} This follows easily from the fact that $$\hom_{K^{-}(B_{0})}(P_{0},I\tensor 
P_{0}[1])\to\ext^{1}(E_{0},E_{0}\ltensor_{A_{0}}I)$$ is an isomorphism 
(where the left-hand side is homotopy classes of maps of the complexes).
\end{proof}

\subsection{The general case}

\begin{lem}\label{L:localness} The condition that a triangle $E_{0}\ltensor_{A_{0}}I\to E\to 
E_{0}\to $ be a deformation of $E_{0}$ is local on $X$.
\end{lem}
\begin{proof} This follows from the fact that formation of the 
natural triangle is natural (!) and the fact that a given map $E\to 
F$ in $\D(X_{A_{0}})$ is an isomorphism if and only if it is 
everywhere locally.
\end{proof}

\begin{proof}[Proof of \ref{T:deformation theory}]  We claim that the 
obstruction $\omega(E_{0})$ is given by the coboundary 
$\delta\theta(E_{0})$.  Indeed, by 
\ref{L:deformations live here}, if there is a deformation, then 
$\delta\theta=0$.  Conversely, if $\delta\theta=0$, then there is 
some triangle $E_{0}\ltensor_{A_{0}}I\to E\to E_{0}\to$ giving rise 
to $\theta$ under the coboundary map.  This remains true affine 
locally on $X$, and we see from \ref{P:affine deformation theory} and 
\ref{L:localness} that $E$ is a deformation of $E_{0}$.  To see that 
the space of deformations is as claimed, the reasoning is similar: one 
knows from \ref{L:deformations live here} that the deformations all 
lie in the fiber over $\theta$.  Reducing to a local computation by 
\ref{L:localness} and using \ref{P:affine deformation theory} shows 
that every element of the fiber is a deformation.  The rest follows 
from \ref{L:splitness}.
\end{proof}

\subsection{Small affine pushouts}

Ultimately, to apply Artin's Representability Theorem, we need to 
prove a version of the Schlessinger-Rim criteria.  In this section we 
will again show that these criteria hold by reducing to the affine 
case and then using the explicit deformation theory provided by 
\ref{C:explicit affine}.

\begin{prop}\label{P:schlessinger} Let 
    $$\xymatrix{ & B\ar[d] \\
    A\ar[r] & A_{0}}$$
    be a diagram of Noetherian $S$-rings with $A\to A_{0}$ a 
    square-zero extension.  Let $E_{0}$ be a complex on $X_{A_{0}}$.  
    Given a deformation $E_{A}$ of $E_{0}$ to $X_{A}$ and a lift 
    $E_{B}$ of $E_{0}$ to $X_{B}$, there is a complex 
    $E\in\D^{b}(X_{A\times_{A_{0}}B})$ such that 
    $E\ltensor_{A\times_{A_{0}}B}A\cong E_{A}$ and $E\ltensor_{A\times_{A_{0}}B}B\cong E_{B}$.
\end{prop}
\begin{proof} The idea of the proof is to make an abstract complex and 
check that it achieves the desired result by reducing to the case of 
affine $X$.  To this end, write $C:=A\times_{A_{0}}B$, and let 
$i_{A_{0}}:X_{A_{0}}\inj X_{C}$, $i_{B}:X_{B}\to X_{C}$, and 
$i_{A}:X_{A}\to X_{C}$ be the natural maps.  We claim that the homotopy 
fiber $F$ of the map $\R(i_{A})_{\ast}E_{A}\oplus
\R(i_{B})_{\ast}E_{B}\to\R(i_{A_{0}})_{\ast}E_{0}$ in $\D(X_{C})$ 
gives the complex we seek.  To prove this, note that there are 
natural maps $F\ltensor_{C}A\to E_{A}$ and $F\ltensor_{C}B\to E_{B}$ 
arising from the adjunction of $\bullet\ltensor_{C}A$ and 
$\R(i_{A})_{\ast}$ 
(and similarly for $B$).  To show that these natural maps are 
quasi-isomorphisms, it suffices to show this locally on $X$; thus, we 
may assume that $X$ is affine.  We may now resolve $E_{B}$ by a 
bounded above complex $P_{B}$ of free $\ms O_{X_{B}}$-modules.  Thus, $E_{0}$ 
is represented by $P_{0}:=P_{B}\tensor_{B}A_{0}$.  Letting 
$P^{i}_{A}$ be the natural lift of $P_{0}^{i}$ to a free 
$\ms O_{X_{A}}$-module, we see from \ref{C:explicit 
affine} that the deformation $E_{A}$ arises from a lift of the 
differentials $\bar d^{i}:P_{0}^{i}\to P_{0}^{i+1}$ to differentials 
$d^{i}:P_{A}^{i}\to P_{A}^{i+1}$, yielding a complex $P_{A}$ such 
that $P_{A}\tensor_{A}A_{0}=P_{0}$.  Furthermore, as affine morphisms 
are cohomologically trivial, we see that the complexes $P_{B}$ and 
$P_{A}$, viewed as $\ms O_{X_{A_{0}}}$-modules by restriction of scalars, 
compute $\R(i_{B})_{\ast}E_{B}$ and $\R(i_{A})_{\ast}E_{A}$.  The natural 
maps $\R(i_{B})_{\ast}E_{B}\to\R(i_{A_{0}})_{\ast}E_{0}$ is realized 
by $P_{B}\to P_{0}$ and similarly for $A$.  We can now apply the 
basic Lemma 3.4 of Schlessinger \cite{schl} to see that the fiber product 
$P_{B}^{i}\times_{P_{0}^{i}}P_{A}^{i}$ is a flat $C$-module.  Thus, 
the complex $P_{B}\times_{P_{0}}P_{A}$, which fits into an exact 
sequence
$$0\to P_{B}\times_{P_{0}}P_{A}\to P_{B}\oplus P_{A}\to P_{0}\to 0,$$
is composed of modules which may be used to compute the derived 
functors $\bullet\ltensor_{C}A$ and $\bullet\ltensor_{C}B$.  The rest 
follows as in Schlessinger to show that 
$(P_{B}\times_{P_{A_{0}}}P_{A})\tensor_{C}A=P_{A}$ and similarly for 
$B$, yielding the result.
\end{proof}

\subsection{Formal deformations}\label{S:groth exist}
In this section, we prove the Grothendieck Existence Theorem for formal deformations over the formal 
spectrum of a complete local Noetherian ring.  While our methods 
will not extend to the arbitrary adic case, what we prove here will 
suffice as input into Artin's theorem in section \ref{S:algebraicity}.

\begin{prop}\label{P:groth exist} Let $(A,\mf m,k)$ be a complete local 
Noetherian $S$-ring and $E_{i}\in\D^{b}(X_{A/\mf m^{i+1}})$ a system of 
elements with compatible isomorphisms $$E_{i}\ltensor_{A/\mf 
m^{i+1}}A/\mf m^{i+2}\simto E_{i+1}$$ in the derived category.  Then 
there is $E\in\D^{b}(X_{A})$ and compatible isomorphisms 
$E\ltensor_{A}A/\mf m^{i+1}\simto E_{i}$.
\end{prop}

Our strategy of proof will (unfortunately) be to express the 
explicit deformation theory in one more way: starting with an 
\emph{injective\/} resolution of $E_{0}$ over $X_{k}$.  Our proof has 
the flaw that it relies heavily on the fact that the ``initial step'' 
takes place over a field, so that $I$ is a free $k$-module in making 
the small extensions, etc.  Once we have done this, we will be able to replace the formal deformation of $E_0$ in the derived category by a formal deformation of complexes.  By reducing to the affine case, we will be able to check that the inverse limit of this deformation of complexes gives a complex on the formal scheme $\widehat X$ which is relatively perfect and whose associated formal deformation (in the derived category) is the original one.  On the other hand, any relatively perfect complex on $\widehat X$ may be constructed by forming finitely many distinguished triangles starting with its (coherent) cohomology sheaves.  Applying the classical form of Grothendieck's Existence Theorem for coherent sheaves will algebraize the inverse limit complex and this will complete the proof.

\begin{lem}\label{L:deformation is deformation} Let $A$ be a local Artinian ring with residue field $k$.  
Suppose $J^{-1}\to J^{0}\to J^{1}$ is an exact sequence of flat 
$A$-modules such that $\H^{0}(J\tensor k)=0$.  Then $\im J^{-1}$ is 
flat and of formation compatible with base change.
\end{lem}
\begin{proof} This resembles \ref{L:non-triv groth lem}, but there 
are no finiteness conditions on the modules $J^{i}$, owing to the 
fact that $A$ is Artinian, as we will see in a moment.  Note that if 
$0\to M\to N\to P\to 0$ is an exact sequence of $A$-modules, then the 
flatness of the $J^{i}$ yields an exact sequence $$0\to J\tensor M\to 
J\tensor N\to J\tensor P\to 0$$ of complexes, yielding an exact 
sequence $$\H^{0}(J\tensor M)\to\H^{0}(J\tensor N)\to\H^{0}(J\tensor 
P).$$  Thus, starting with the vanishing of $\H^{0}(J\tensor k)$ and 
proceeding by induction on length, we conclude that $\H^{0}(J\tensor 
M)=0$ for any finite $A$-module $M$ (and then for any $A$-module by 
direct limit considerations).  In particular, $\H^{0}(J)=0$, yielding 
an exact sequence
$$0\to J^{0}/\im J^{-1}\to J^{1}\to\coker(J^{0}\to J^{1})\to 0$$
the universal vanishing of $\H^{0}(J\tensor M)$ shows that this 
sequence remains exact upon tensoring with any $A$-module $M$, whence 
we conclude that all terms are flat.  It follows that $\im 
J^{-1}$ is flat over $A$ and that $(\im J^{-1})\tensor M=\im (J^{-1}\tensor M)\subset 
J^{0}\tensor M$.
\end{proof}

\begin{lem}\label{L:local computation of pullback} Let $f:(X,\ms 
O_{X})\to (Y,\ms O_{Y})$ be a morphism of ringed 
topoi.  Let $J$ be a complex of $\ms O_{Y}$-modules.  Suppose there 
exists a covering $U_{i}\to e_{Y}$ of the final object of $Y$ such 
that for each $i$, there is a $K$-flat complex $\ms J_{i}$ on 
$U_{i}\in Y$ and a map $\psi:\ms J_{i}\to J|_{U_{i}}$ such 
that both $\psi$ and $f^{\ast}\psi:f^{\ast}\ms J_{i}\to 
f^{\ast}J|_{U_{i}}$ 
are quasi-isomorphisms.  Then the complex $f^{\ast}J$ computes the 
derived pullback ${\field L} f^{\ast}J$.
\end{lem}
\begin{proof} Let $F\to J$ be a termwise surjective $K$-flat 
resolution (constructed e.g.\ using \ref{L:natural good}).  It is easy (from 
the existence and exactness of extension by zero from $U_{i}$ to $X$) to 
see that $F|_{U_{i}}\to J|_{U_{i}}$ is also a $K$-flat 
resolution, and that there exists a commutative diagram
$$\xymatrix{& F|_{U_{i}}\ar[dr] & \\ 
G\ar[ur]\ar[dr] & & J|_{U_{i}}\\
&  \ms J_{i}\ar[ur] &}$$
of quasi-isomorphisms with all but the right-most complex $K$-flat.  
By the general theory of $K$-flatness, the two left-hand arrows pull 
back to give quasi-isomorphisms on $X$.  By hypothesis, the bottom 
right arrow also pulls back to give a quasi-isomorphism.  The result 
follows from the fact that pulling back $K$-flat complexes computes 
the derived pullback.
\end{proof}

Let $E_{0}\to J$ be a resolution of $E_{0}$ by a bounded below complex 
of injective $X_{0}$-modules.  

\begin{prop}\label{P:injective resolution} Let $(A_{0},\mf m,k)$ be a local Artinian $S$-ring, $A\to A_{0}$ a 
local small extension with kernel $I$ annihilated by $\mf m$.  Let 
$J_{0}$ 
be a bounded below complex of $A_{0}$-flat $X_{A_{0}}$-modules such 
that  
$J_{0}\tensor_{A_{0}}k$ is a complex of injective $X_{k}$-modules 
which computes $J_{0}\ltensor_{A_{0}}k$.  Then any 
deformation of $J_{0}$ to $X_{A}$ has the form $J$, with $J^{i}$ an 
$A$-flat deformation of $J_{0}^{i}$ for each term $i$.  In 
particular, one has that the complex $J_{0}$ computes 
$J\ltensor_{A}A_{0}$.  Furthermore, the 
termwise deformations $J_{0}^{i}$ are unique up to isomorphism.
\end{prop}
\begin{proof} Since $J_{k}^{i}$ is injective and $I$ is a free 
$k$-module, the usual deformation theory for modules shows everything 
about the termwise deformations: there is an 
obstruction in $\ext^{2}(J_{k}^{i},I\tensor J_{k}^{i})$ and deformations are 
parametrized by $\ext^{1}(J_{k}^{i},I\tensor J_{k}^{i})$, both of 
which vanish (as $J_{k}^{i}\tensor I$ is still injective by freeness 
of $I$ over $k$).  Thus, we may let $J^{i}$ denote the essentially unique 
$A$-flat lift of $J_{0}^{i}$ to $X_{A}$.  A similar computation shows 
that the differential $\bar d^{i}:J_{0}^{i}\to J_{0}^{i+1}$ admits a 
lift $d^{i}:J^{i}\to J^{i+1}$.
Just as in \ref{Constr:obstruction}, choosing lifts yields an 
obstruction in $\ext^{2}(J_{k},I\tensor J_{k})$.  

Using \ref{R:from projective to flat}, it is easy to see that the 
obstruction just constructed equals $\omega(E_{0})$.  Indeed, one may 
assume that the lifts $F^{i}\to F^{i+1}$ of the 
differentials $F_{0}^{i}\to F_{0}^{i+1}$ may be chosen to cover a 
lift $J^{i}\to J^{i+1}$ (for a suitable choice of good resolution 
$F_{0}\to J_{0}$).  To see this, first note that $J^{i}\to 
J^{i}_{0}$ is surjective on sections over any affine (as $I$ is a 
free $k$-module, so $J^{i}_{0}\tensor I$ is an injective sheaf on 
$X\tensor k$, hence has vanishing $\H^{1}$).  Thus, if we take the 
natural lift $F^{i}$ of the component $F^{i}_{0}$ of the canonical 
good resolution $F_{0}\to J_{0}$, we can lift the map $F_{0}^{i}\to 
J_{0}^{i}$ in various ways to a map $F^{i}\to J^{i}$.  
Let $F_{0}\to J_{0}$ be the canonical good resolution of \ref{L:natural 
good}; in particular, we may assume that $F_{0}\to J_{0}$ is surjective on 
sections over any \'etale affine $V\to X$.  It follows that 
$F_{0}\tensor k\to J_{0}\tensor k$ is surjective on sections over any 
affine, and thus that $I\tensor F_{0}\to I\tensor J_{0}$ is 
surjective on sections over any affine, as $I$ is a free $k$-module.  The 
chosen map $F^{i}\to J^{i}$ is 
easily seen to be surjective by a simple Snake Lemma argument.  
Moreover, it is easy to see that the map $F^{i}\to 
F^{i}_{0}$ is surjective on sections over any affine (as a section 
of $\oplus j^{V}_{!}\ms O_{V}$ over $U$ is given by a factorization 
$U\to V$ over $X$ and a section of $\ms O_{U}$).  We claim that the induced map $F^{i}\to 
J^{i}\times_{J_{0}^{i}}F_{0}^{i}$ is surjective on sections over any 
affine.  This is easily deduced from the diagram
$$\xymatrix{0\ar[r] & I\tensor F_{0}^{i}(V)\ar[r]\ar[d] & F^{i}(V)\ar[r]\ar[d] & 
F_{0}^{i}(V)\ar[r]\ar[d] & 0\\
0\ar[r] & I\tensor J_{0}^{i}(V)\ar[r] & J^{i}(V)\ar[r] & 
J_{0}^{i}(V)\ar[r] & 0.}$$
From this it follows that the map $F_{0}^{i-1}\to 
F_{0}^{i}$ may be lifted to a map $F^{i-1}\to F^{i}$ over the chosen 
lift $J^{i-1}\to J^{i}$, as the lift is given by a collection of 
maps $j^{V}_{!}\ms O_{V}\to F^{i}$, which amounts to a collection of 
sections of $F^{i}$ over $V$.

Thus, if $E_{0}$ is unobstructed, there also exists a lift of $J_{0}$ 
to a complex $J$; furthermore, the space of such lifts is a torsor 
under $$\ext^{1}_{X_{A_{0}}}(J_{0},I\tensor 
J_{0})=\ext^{1}_{X_{k}}(J_{0}\tensor k,I\tensor J_{0}\tensor 
k)=\ext^{1}_{X_{k}}(E_{0}\ltensor k,I\ltensor E_{0}\ltensor k).$$
If we can show that $J\ltensor_{A}A_{0}\to J_{0}$ is a 
quasi-isomorphism, then we will have shown that all deformations arise 
as lifts to complexes $J$.  To see this last point, we use the 
boundedness of $J_{0}$ and \ref{L:deformation is deformation} to 
conclude that for sufficiently large $n$, one has $\tau_{\leq n}J\to J$ is a quasi-isomorphism of complexes of flat modules over $A$ and $\tau_{\leq 
n}J\tensor_{A}A_{0}=\tau_{\leq n}J_{0}$ as complexes of flat modules 
over $A_{0}$.  Since bounded complexes of flat 
modules compute the derived tensor product, it immediately follows 
that $J_{0}$ computes $J\ltensor_{A}A_{0}$ and thus that $J$ is a 
deformation of $J_{0}$.
\end{proof}

\begin{cor}\label{C:formal complex} In the situation of \ref{P:groth exist}, the formal 
deformation $(E_{i})$ is the image of a formal deformation $J_{i}$ of 
quasi-coherent complexes.
\end{cor}

\begin{remark}\label{R:local computation} Note that combining \ref{P:injective resolution} with 
\ref{C:explicit affine} easily yields further information: there is an 
affine covering $\widehat U_{i}$ of the formal completion $\widehat X$ of 
$X$ along the closed fiber such that there is a bounded above complex 
of finite free $\ms O_{\widehat X}$-modules $P$ and a system of quasi-isomorphisms 
$P_{i}\to(J_{i})|_{\widehat U}$.  (This will be taken up again in 
\ref{L:local coho} below for the reader desiring clarification.)  We will see in a moment that this 
implies the hypothesis of \ref{L:local computation of pullback} for 
the system of sheaves $\invlim J_{i}$ on $\widehat X$.  This adds a 
certain amount of naturality to some of the following constructions, 
in the sense that rather than constantly paying attention to a truncation 
as in the proof of \ref{P:injective resolution}, we can simply work 
with the system of unbounded injective resolutions $(J_{i})$ without 
fear of losing control over the derived pullback to the various 
thickenings of $X$.
\end{remark}

\begin{lem}\label{L:completion lemma} Let $A$ be a Noetherian ring, $\mf m\subset A$ an ideal,  
and $(M_{i})$ a system of modules 
over $(A_{i}=A/\mf m^{i+1})$ such that $M_{i}\tensor A_{j}=M_{j}$ for 
all $j\leq i$.  Let $M=\invlim M_{i}$.  The natural map $A_i\tensor M\to M_i$ is an isomorphism.
\end{lem}
\begin{proof} Clearly the map $M\tensor A_{i}\to M_{i}$ is 
surjective.  To show that it is an isomorphism, it is thus enough to 
prove that any element of $M$ which maps to zero in $M_{i}$ is in 
$\mf m^{i+1}M$.  Since $A$ is Noetherian, $\mf m^{i+1}$ is finitely 
generated, say by $v_{0},\ldots,v_{n}$.  An element $m$ in the kernel of 
$M\to M_{i}$ corresponds to a system of elements $m_{j}\in M_{j}$ 
such that $m_{j}=0$ for all $j\leq i$ and $\widebar{m_{j+1}}=m_{j}$, 
where the bar denotes the reduction map $M_{j+1}\tensor A_{j}\simto 
M_{j}$.  Thus, we conclude first of all that for all $j$, 
$m_{j}\in\mf m^{i+1}M_{j}$.  We wish to write the element $m$ in the 
form $\sum v_{\alpha}m^{\alpha}$ for $m^{\alpha}\in M$.  We can 
certainly do this in $M_{i+1}$ by the assumptions about the system 
$(M_{i})$.  Suppose we have done this compatibly for $m_{k}$ for all $k\leq j$.  
Since $M_{j+1}\to M_{j}$ is surjective, we may lift 
each $m_{j}^{\alpha}$ to some element $\tilde m_{j}^{\alpha}\in 
M_{j+1}$.  In this case, we have that $b:=m_{j+1}-\sum v_{\alpha}\tilde 
m_{j}^{\alpha}$ maps to $0\in M_{j}$.  But the kernel of 
$M_{j+1}\to M_{j}$ is $\mf m^{j+1}M_{j+1}$, so the kernel of $\mf 
m^{i+1}M_{j+1}\to\mf m^{i+1}M_{j}$ is $\mf m^{j-i}\mf m^{i+1}M_{j+1}$.  Thus, we may write $b=\sum v_{\alpha}\beta^{\alpha}$, where $\beta^{\alpha}\in\mf 
m^{j-i}M_{j+1}$.  Setting $m^{\alpha}_{j+1}=\beta^{\alpha}+\tilde m_{j}^{\alpha}$, we see that $m^{\alpha}_{j+1}$ agrees with $m^{\alpha}_j$ in $M_{j-i}$.  Thus, the terms stabilize, and by induction we are done.
\end{proof}

\begin{lem} Let $(A,\mf m)$ be as above and $X\to A$ a Noetherian algebraic space, with $X_i=X\tensor_A A_i$.  Let $(J_{i})$ be a system of quasi-coherent 
sheaves on $X_{i}$ such that for all $i\leq j$ one has $J_{i}\tensor 
A_{j}=J_{j}$.  Then the inverse limit $\invlim J_{i}$ defines a (not 
necessarily quasi-coherent!) sheaf 
of $\ms O_{\widehat X}$-modules $J$ such that $J\tensor_{\ms O_{\widehat 
X}}\ms O_{X_{i}}=J_{i}$ for all $i$.
\end{lem}
\begin{proof} By \ref{L:completion lemma}, we see that for any \'etale affine 
$U\to X$, the natural map $$J(\widehat U)\tensor_{A}A_{i}\to 
J(U_{i})$$ is an isomorphism.  Since affines of the form $\widehat U$ form a basis for the 
\'etale topology on $\widehat X$, the result follows by sheafification.
\end{proof}

\begin{defn} A system $(J_{i})$ of quasi-coherent sheaves on $X_{i}$ 
such that $J_{i}|_{X_{j}}=J_{j}$ for all $j\leq i$ will be called an 
\emph{ind-quasi-coherent\/} sheaf on $\widehat X$.  If each sheaf $J_i$ is flat over $A_i$, then we will say that $(J_i)$ is \emph{ind-flat\/} over $\spf A$.
\end{defn}

It is clear by definition that the mapping cone of a map of complexes 
of (ind-flat) ind-quasi-coherent sheaves on $\widehat X$ is itself 
(ind-flat) ind-quasi-coherent.  The local criterion of flatness shows that any ind-flat sheaf on $\widehat X$ which is of finite presentation is also flat over $\spf A$; unfortunately, it is rarely the case that ind-quasi-coherent sheaves are of finite presentation, so we cannot assume that ind-flat sheaves are flat.  This will turn out not to make a difference for us.

\begin{lem}\label{L:formal comparison} Let $X\to S=\spec R$ be a scheme over a complete local Noetherian 
ring.  Let $F^{\bullet}$ and $J^{\bullet}$ be complexes of 
$\spf R$-ind-flat
ind-quasi-coherent sheaves on the formal completion $\widehat X$ along the closed 
fiber.  If a map $\phi:F^{\bullet}\to J^{\bullet}$ of complexes has 
the property  
that for all $n$, $\phi_{n}:F^{\bullet}_{n}\to J^{\bullet}_{n}$ is a 
quasi-isomorphism, then $\phi$ is quasi-isomorphism.
\end{lem}
\begin{proof} Taking the mapping cone of $\phi$ and noting that 
its formation commutes with reduction modulo powers of $\mf m_{R}$, it 
suffices to prove that if $K^{\bullet}$ is a complex of ind-flat  
ind-quasi-coherent $\ms O_{\widehat X}$-modules such that $K^{\bullet}_{n}$ is 
acyclic, then $K^{\bullet}$ is acyclic.  By construction we have that 
the induced maps $K_{n+1}\to K_{n}$ are termwise surjective with 
kernel $I_{n+1}\tensor K_{0}$.  It follows by an easy diagram chase 
that $K_{n+1}$ is acyclic, and then that $K$ is acyclic on $\widehat X$.
\end{proof}

\begin{lem}\label{L:local coho}  Let $X\to S=\spec R$ be a flat map 
from an affine scheme to a complete local Noetherian ring with 
completion $\widehat X$ along the closed fiber.  Suppose 
$J^{\bullet}$ is a complex of $\spf(R)$-ind-flat ind-quasi-coherent $\ms O_{\widehat 
X}$-modules such that for all $i$, $J^{\bullet}\tensor R/\mf m^{i}=J^{\bullet}_{i}$ is relatively perfect and bounded as an object of $\D(X\tensor R/\mf m^{i})$.  Then $J^{\bullet}$ is bounded with coherent cohomology and $J^{\bullet}\ltensor R/\mf m^i\cong J^{\bullet}_i$ in $\D(X\tensor R/\mf m_i)$.
\end{lem}
\begin{proof} We may choose a resolution $\phi_{0}:F^{\bullet}_{0}\to 
J_{0}^{\bullet}$ consisting of a bounded above complex of finite free 
$\ms O_{X_{0}}$-modules.  By \ref{C:explicit affine}, there will be a 
lift $F_{0}$ to $F_{1}$ and a quasi-isomorphism $\phi_{1}:F_{1}\to J_{1}$ lifting 
$F_{0}\to J_{0}$ up to a homotopy.  Since the components of $F_{1}$ 
are free, we may clearly lift the homotopy and conclude that 
$\phi_{1}$ is a lift of $\phi_{0}$ in the category of complexes (not 
up to homotopy).  Continuing in this manner yields a map $\phi:F\to J$ on 
$\widehat X$, where $F$ is a bounded above complex of finite free $\ms O_{\widehat 
X}$-modules, such that $\phi_{i}$ is a quasi-isomorphism for all 
$i$.  Applying \ref{L:formal comparison} and \ref{L:local computation of pullback} finishes the proof (once we 
note that the deformations in $\ppD$ cannot grow new cohomology 
sheaves).
\end{proof}

\begin{proof}[Proof of \ref{P:groth exist}] By \ref{L:local computation of pullback}, \ref{P:injective 
resolution}, and \ref{L:local coho}, we have that there exists a complex $E$ of ind-flat ind-quasi-coherent modules on $\widehat X$ and isomorphisms $\mathbf{L}\iota_i^{\ast}E\simto E_i$, where $\iota_{i}:X\tensor A/\mf m^{i+1}\inj\widehat X$ is the natural closed immersion of formal schemes.  Furthermore, this formal complex $E$ has bounded coherent cohomology.  By 
the classical Grothendieck Existence Theorem (\S 5 of \cite{ega3-1}), each of the 
cohomology sheaves is algebraizable, and for any coherent sheaves $F$ 
and $G$ on $X$, one has $\ext_{X}^{i}(F,G)=\ext^{i}_{\widehat X}(\widehat 
F,\widehat G)$.  It follows by induction that $E$ is algebraizable.  (The 
reader will note that the map $\widehat X\to X$ is faithfully flat on the category of 
coherent sheaves, hence the ``formalization'' of an algebraic 
coherent sheaf is equal to its left-derived formalization.  We are 
implicitly using the fact -- Corollaire 2.2.2.1 of Expos\'e II of \cite{sga6} -- that 
$\D^{b}(\coh(X))=\D^{b}_{\coh}(X)$, and similarly for $\widehat X$.)
\end{proof}

\section{Algebraicity}\label{S:algebraicity}

Having assembled the necessary preliminaries, we are now ready to 
prove that the stack of universally gluable relatively perfect 
complexes is algebraic in the sense of Artin.

\subsection{Constructibility properties}

Here we verify that the deformation 
theory and obstruction theory for relatively perfect complexes are 
constructible in the sense of 4.1 of \cite{artin}.

\begin{lem}\label{L:constructibility} Let $X\to\spec A_{0}$ be a 
proper flat algebraic space of finite 
presentation over a reduced Noetherian ring and $I$ a finite 
$A_{0}$-module.  Let $A_{0}\to B_{0}$ be an \'etale ring extension.  Given $E\in\ppD(X/A_{0})$ and any $i$, 
\begin{enumerate}
    \item $\ext^{i}(E,E\ltensor I)\tensor_{A_{0}} 
    B_{0}=\ext^{i}(X_{B_{0}}(E_{B_{0}},I_{B_{0}}\ltensor E_{B_{0}})$;

    \item if $\mf m\subset A_{0}$ is a maximal ideal then 
    $$\ext^{i}_{X}(E,E\ltensor I)\tensor_{A_{0}}\widehat 
    A_{0}=\invlim\ext^{i}_{X}(E,I/\mf m^{i+1}I\ltensor E);$$

    \item if $X$ is a scheme then there is a dense set of 
points $p\in\spec A_{0}$ such that the natural map 
$$\ext^{i}(E,E\ltensor I)\tensor\kappa(p)\to\ext^{i}_{X_{\kappa(p)}}(E_{p},I_{p}\ltensor_{X_{\kappa(p)}}E_{p})$$
is an isomorphism.
\end{enumerate}
\end{lem}
\begin{proof} Part (1) is trivial.  
    
To prove (2), we may replace $A_{0}$ with $\widehat{(A_{0})}_{\mf m}$ and assume that $A_{0}$ is a 
complete local Noetherian ring.  We claim that for any $i$, the natural map 
$$\rho:\sext^{i}(E,I\ltensor E)\to\invlim_{j}\sext^{i}_{X_{j}}(E_{j},I/\mf 
m^{j+1}\ltensor E_{j}),$$
where the subscript $j$ denotes the derived base change to $A_{0}/\mf 
m^{j+1}$.  Since $\rho$ is natural in $X$, to prove that it is an 
isomorphism we may assume that $X$ is affine.  In this case, we can 
choose two resolutions of $E$: one $Q^{\bullet}\to E$ by a bounded 
complex of $A_{0}$-flat coherent sheaves and one $P^{\bullet}\to E$ by 
a bounded above complex of free $\ms O_{X}$-modules of finite rank.  Then 
the complex $\ms C:=\shom(P^{\bullet},Q^{\bullet}\tensor I)$ computes the 
left hand-side of $\rho$ and it is easy to see that $\ms 
C\tensor_{A_{0}}A/\mf m^{j+1}$ computes the $j$th term of the right-hand 
side.  
Since $\ms C$ is a complex of 
coherent $\ms O_{X}$-modules, the result follows by standard 
Mittag-Leffler arguments.  We claim that once we know $\rho$ is an 
isomorphism, the proof of (2) follows immediately.  Indeed, by the 
Grothendieck comparison theorem \cite{ega3-1} and the result just 
proved, we have that $\H^{p}(\sext^{q}(E,E\ltensor 
I))=\invlim_{j}\H^{p}(\sext^{q}(E_{j},E_{i}\ltensor I/\mf 
m^{j+1}I))$.  On the other hand, there are the usual spectral sequences 
$E^{pq}_{2}=\H^{p}(\sext^{q}(E,E\ltensor I))\Rightarrow\ext^{p+q}(E,E\ltensor 
I)$ and $E^{pq}_{j,2}=\H^{p}(\sext^{q}(E_{j},E_{j}\ltensor I/\mf 
m^{j+1}I))\Rightarrow \ext^{p+q}(E_{j},E_{j}\ltensor I/\mf 
m^{j+1}I)$, and these are compatible with the maps in the system 
$\rho$.  Thus, to prove (2) it is enough to prove that the abutment of 
the inverse limit of the spectral sequences $E^{pq}_{j,\bullet}$ is 
the inverse limit of the abutments.  This follows from a standard 
Mittag-Leffler argument (made extremely easy by the fact that the 
$E^{pq}_{j,r}$ vanish for sufficiently large $p$ independent of $q$, 
$j$, and $r$).

To prove (3), one reasons just as in \ref{P:locally quasi-sep 
essence}: one first shrinks $S$ until $I$ is flat, then descends to a regular finite-type $\Z$-algebra and 
then finds $S^{0}\subset S$ over which $\tau_{\leq n}\R 
f_{\ast}\rshom(E,E\ltensor I)$ universally computes (upon any derived 
pullback) the desired $\ext$ module (for some sufficiently large $n$).  
Shrinking $S^{0}$ further so that all of the cohomology sheaves of 
$\tau_{\leq n}\R 
f_{\ast}\rshom(E,I\ltensor E)$ (which is now bounded) are $S$-flat yields an open set 
satisfying the required conditions.  We have used the fact that $X$ is 
a scheme only in descending the complex $E$ to a finite 
$\Z$-subalgebra of $A_{0}$, where we have to invoke \ref{P:local finite presentation}, whose 
proof (unfortunately) used schemehood.
\end{proof}

\subsection{Proof of algebraicity}

\begin{thm}\label{T:main theorem} Let $f:X\to S$ be a flat proper morphism of finite 
presentation between algebraic spaces which is fppf-locally on $S$ representable by 
schemes.  The stack $\sppugD(X/S)$ is an Artin stack locally of 
finite presentation (and locally quasi-separated) over $S$.
\end{thm}
This seems to be a basic object which will be useful for the future study 
of many moduli problems arising from perverse coherent sheaves and other 
natural structures on derived categories.  In particular, one hopes to 
be able to better understand Bridgeland's construction \cite{bridgeland} 
using the general theory of this paper.  It also seems as though this 
theory should be useful for understanding the results of Abramovich and 
Polishchuk on valuative criteria for stable complexes 
\cite{abramovich-polishchuk}.

\begin{proof} The proof consists of checking a few conditions, which 
we have painstakingly proven above!  The local finite presentation 
condition is \ref{P:local finite presentation}.  The Schlessinger 
conditions are \ref{P:schlessinger}.  The deformation theory of 
objects and automorphisms is 
described by \ref{T:deformation theory} (with obvious linearity and 
functoriality), and its contructibility properties are 
\ref{L:constructibility}.  Finally, the local quasi-separation is 
\ref{C:local quasi-sep}.  In fact, as is shown in \ref{P:locally 
quasi-sep essence}(2), the diagonal is actually separated, which is usually required 
for the stack to be considered algebraic in modern terminology.
\end{proof}


\subsection{Easy applications}\label{S:easy applications}

Using standard results about ``rigidifications'' of algebraic stacks (dating back to item 
2 of the appendix to Artin's \cite{artin}  
and rediscovered by Abramovich, Corti, and Vistoli in 
\cite{abramovich-corti-vistoli}), we may produce algebraic spaces representing sheafified 
functors of complexes with only scalar automorphisms.  This will 
reproduce the main result of Inaba's paper \cite{inaba} (and extend it to the 
case of an arbitrary flat proper representable morphism of finite presentation 
between quasi-separated algebraic spaces, rather than simply a 
flat projective morphism of Noetherian schemes).

\begin{defn} A complex $E$ is \emph{simple\/} if $E\in\ppugD(X/S)$ 
and the natural map $\G_{m}\to\saut(E)$ is an isomorphism.
\end{defn}

It is clear that simple universally gluable $S$-perfect 
complexes form a stack $s\!\sppugD(X/S)$ with a natural substack 
structure $\eta:s\!\sppugD\to\sppugD$.  

\begin{lem} With the above notation, the map $\eta$ is an open 
immersion.  
\end{lem}
\begin{proof} Since $\sppugD$ is an algebraic stack, we know that the inertia stack $\ms I\to\sppugD$ is represented by finite type separated morphisms of schemes.  Moreover, the fibers are all naturally open subsets of affine spaces, hence are equidimensional.  It is easy to see that $s\!\sppugD$ is identified with the open substack corresponding to the locus where the fibers of $\ms I\to\sppugD$ have dimension at most $1$ at any point.
\end{proof}

\begin{cor} The stack $s\!\sppugD(X/S)$ is an Artin stack locally of 
finite presentation over $S$.  There is an algebraic space locally of 
finite presentation $D^{s}$ and a morphism $s\!\sppugD\to D^{s}$ 
which is a $\G_{m}$-gerbe.  In particular, $D^{s}$ represents the 
sheafification of $s\!\sppugD$ on the big \'etale site of $S$.
\end{cor}
\begin{proof} It remains to find the space $D^{s}$.  Since the inertia 
stack of $s\!\sppugD$ is naturally identified with $\G_{m}$, this 
follows immediately from the theory of rigidification of 
Abramovich-Corti-Vistoli \cite{abramovich-corti-vistoli}, or (in the case of simple 
sheaves) from an argument of Artin (contained in the appendix to \cite{artin}).
\end{proof}

\begin{cor}[Inaba] Let $X\to S$ be a projective morphism of Noetherian 
schemes.  Define a functor $F$ on $S$-schemes as follows: to $T\to S$, 
associate the set of quasi-isomorphism classes of simple bounded complexes 
of $S$-flat coherent sheaves on $X$.  Then the sheafification of $F$ 
is representable by an algebraic space locally of finite type over $S$.
\end{cor}
\begin{proof} The careful reader will invoke the following result 
(II.2.2.2.1 of \cite{sga6}): $$\D^{b}_{\coh}(X)=\D^{b}(\coh(X)).$$
\end{proof}

\end{document}